\documentstyle[epsf,stmaryrd,amscd,amstex,amssymb,
11pt
]{amsart}
\input xy

\newcommand{\cZ}{{\cal  Z}}

\newcommand{\fU}{{\mathfrak U}}

\newcommand{\ot}{\otimes}

\newcommand{\<}{\langle}

\renewcommand{\>}{\rangle}

\newcommand{\emp}{\emptyset}

\newcommand{\fZ}{{\mathfrak Z}}

\newcommand{\HH}{{\cal H}}

\newcommand{\VV}{{\cal V}}

\newcommand{\SS}{{\cal S}}
\newcommand{\MM}{{\cal M}}

\newcommand{\WW}{{\mathsf W}}

\newcommand{\cupl}{{\bigcup\limits}}

\newcommand{\isor}{{\widetilde \longrightarrow}}
\newcommand{\st}{\star}

\newcommand{\oplusl}{\bigoplus\limits}

\newcommand{\Ga}{\Gamma}
\newcommand{\Gab}{{\overline{\Gamma}}}
\newcommand{\Gb}{{\overline{G}}}
\newcommand{\Yb}{{\overline{Y}}}
\newcommand{\Tb}{{\overline{\mT}}}
\newcommand{\mTb}{{\overline{\mT}}}

\newcommand{\TO}{{\overline{T}}}

\newcommand{\Zet}{{\mathbb{Z}}}
\newcommand{\Ce}{{\mathbb{C}}}

\newcommand{\La}{{\Lambda}}

\newcommand{\Ohat}{{\widehat{\mathcal O}}}

\newcommand{\Aone}{{\mathbb{A}}^1}

\newcommand{\tN}{\widetilde{N}}

\newcommand{\Nt}{\widetilde{N}}

\newcommand{\Gm}{{\mathbb{G}}_m}

\newcommand{\de}{\partial}

\newcommand{\bbP}{{\mathbb{P}} }

\newcommand{\bbPt}{\tilde{\mathbb{P}} }

\renewcommand{\O}{{\cal O}}

\newcommand{\supp}{\mathop{\mathrm{supp}}}

\newcommand{\A}{{\sf A}} 
\newcommand{\B}{{\sf B}} 
\newcommand{\Be}{{\sf B}}

\newcommand{\fS}{{\mathfrak S}}

\renewcommand{\H}{{\mathcal{H}}}	

\newcommand{\To}{\longrightarrow}

\newcommand{\imbed}{\hookrightarrow}

\newcommand{\la}{\lambda}

\newcommand{\bs}{\backslash}

\newcommand{\Pone}{{\mathbb P}^1}

\newcommand{\Qu}{{\mathbb Q}}

\newtheorem{Thm}{Theorem}
\newtheorem{Cor}[Thm]{Corollary}
\newtheorem{Lem}[Thm]{Lemma}
\newtheorem{Prop}[Thm]{Proposition}
\newtheorem{Claim}[Thm]{Claim}

\theoremstyle{definition}
\newtheorem{Def}[Thm]{Definition}

\theoremstyle{remark}
\newtheorem{Rem}[Thm]{Remark}
\newtheorem*{Rem*}{Remark}
\newtheorem{Ex}[Thm]{Example}

\numberwithin{Thm}{section}

\newcommand{\proof}{{\em Proof \ }}

\def\proof{\smallskip\noindent {\bf Proof.\ }}

\def\BB{{\mathcal{B}}}

\newcommand{\OO}{{\mathcal O}}

\newcommand{\LL}{{\mathcal L}}

\newcommand{\EE}{{\mathcal E}}

\newcommand\nc{\newcommand}
\nc\on{\operatorname}
\nc\renc{\renewcommand}

\nc{\ppart}{(\!(t)\!)}

\nc\Bun{\on{Bun}}
\renc\mod{\text{-mod}}

\nc\CK{{\mathcal K}}
\nc\CO{{\mathcal O}}

\nc\BC{\mathbb C}
\nc\BA{\mathbb A}
\nc\BO{\mathbb O}

\nc{\fg}{{\mathfrak g}}

\nc{\hg}{\widehat\fg}

\newcommand{\mT}{{\mathbb T}}

\newcommand{\mC}{{\mathbb C}}
\newcommand{\mZ}{{\mathbb Z}}
\newcommand{\mP}{\mathbb P}
\newcommand{\mG}{{\mathbb G}}
\newcommand{\mA}{{\mathbb A}}

\newcommand{\GG}{\Gamma}

\newcommand{\bo}{\omega}
\newcommand{\bl}{\lambda}

\newcommand{\bS}{\Sigma}

\newcommand{\mcB}{\mathcal B}

\newcommand{\mcH}{{\mathcal H}}

\newcommand{\mcO}{\mathcal O}

\newcommand{\mcV}{\mathcal V}

\newcommand{\fb}{{\mathfrak b}}

\newcommand{\fl}{\mathfrak l}

\newcommand{\fp}{{\mathfrak p}}

\newcommand{\ft}{\mathfrak t}

\newcommand{\ti}{\tilde}
\newcommand{\un}{\underline}
\newcommand{\bc}{\backslash}

\numberwithin{equation}{section}

\newtheorem{theorem}{Theorem}[section]

\newtheorem{corollary}[theorem]{Corollary}
\newtheorem{lemma}[theorem]{Lemma}
\newtheorem{claim}[theorem]{Claim}

\theoremstyle{definition}

\newtheorem{example}[theorem]{Example}
\newtheorem{definition}[theorem]{Definition}

\newcommand{\fC}{{\mathfrak C}}

\newcommand{\Xb}{\bar{X}}

\newcommand{\Cbar}{\bar{\fC}}

\newcommand{\al}{\alpha}

\newcommand{\aff}{\mathrm{aff}}

\newcommand{\as}{\mathrm{as}}

\begin{document}

\title[Geometry of second adjointness]{Geometry of second adjointness for $p$-adic groups}
\author{Roman Bezrukavnikov, David Kazhdan}

\begin{abstract}
We present a geometric proof of second adjointness for a reductive $p$-adic group. Our approach is based
on geometry of the wonderful compactification and related varieties. Considering asymptotic behavior of
a function on the group in a neighborhood of a boundary stratum of the compactification, we get a "co-specialization" map
between spaces of functions on various varieties carrying a $G\times G$ action. These maps can be viewed as maps of bimodules for the Hecke algebra, and the corresponding natural transformations of endo-functors
of the module category lead to the second adjointness. 
We also get a formula for the "co-specialization" map expressing
it as a composition of the orishperic transform and inverse intertwining operator; a parallel result for $D$-modules was obtained in \cite{BFO}.
 As a byproduct we obtain a formula for the Plancherel functional restricted to a certain commutative subalgebra in the Hecke algebra generalizing a result by Opdam.
\end{abstract}

\maketitle

\centerline{\em To the memory of Izrail' Moiseevich Gel'fand}

\tableofcontents

\section{Introduction}
Parabolic induction and restriction (Jacquet) functors play a fundamental role in representation theory
of reductive $p$-adic groups. 
It follows directly from definitions that the parabolic induction functor
 is right adjoint to the Jacquet functor; we will call this adjointness the
ordinary or 
Frobenius adjointness.
It has been discovered by Casselman 
 for admissible representations and generalized by Bernstein
to arbitrary smooth representations 
 that there is another non-obvious adjointness
 between the two functors. Namely, the parabolic induction functor turns out
 to be also {\em left} adjoint to
Jacquet functor with respect to the opposite parabolic  (we will refer to this as the second or Bernstein adjointness). This fact appears in unpublished
notes of Bernstein \cite{Bern} (see also
exposition in \cite{Rena}), we reprove it below. 
Rather than following the original strategy, our
approach emphasizes the relation to geometry of the group and related spaces.
 More precisely,
showing that two functors $\alpha$ and $\beta$ are adjoint amounts to providing morphisms
between the identity functor and the compositions $\alpha \circ \beta$, $\beta\circ \alpha$ satisfying
certain compatibilities. When $\alpha$, $\beta$ are Jacquet and parabolic induction functors, the arrow between the endo-functors of
the category of representations of the Levi is easy to define, both for the ordinary and for the second adjointness. To describe the morphisms between
the endo-functors of representations of $G$, recall that an endo-functor of a category of modules
usually comes from a bimodule. This is so in our case, moreover, the bimodules can be realized as spaces
of functions\footnote{Here and elsewhere in the Introduction  we ignore the difference between functions and measures, see the  precise statement in section \ref{2ndadj}.} on a variety equipped with a $G\times G$-action. Thus showing adjointness amounts to providing a map
between the function spaces satisfying certain properties.
For the ordinary induction, this turns out to be the map from functions on the group to integral kernels acting on the universal principal series, i.e. on the space of functions on $G/U$. This well studied map (known as the orispheric transform) is given by an explicit correspondence.

In contrast, the map responsible for the second adjointness has not, to our knowledge, received much attention in the literature.\footnote{See, however, \cite{SV} (which appeared after the preprint version of the present paper) where related ideas have been developed and applied in a more general context.} The goal of this paper is a geometric description of this map and its basic properties.
It is a map (to be denoted by $\Be $) from functions on (the space of $F$-points of) a certain algebraic variety $X$ to the space of functions on $G$.
Our first observation is that $X$ is an asymptotic cone of $G$. By this we mean that it admits an open embedding to the normal space $N_{\Gb}(Z)$ where $\Gb$ is the De Concini-Procesi compactification
of $G$ and  $Z\subset \Gb$ is a $G\times G$ orbit.\footnote{When the group $G$ is not of adjoint type
the compactification is not smooth. In this case we use a modification of the notion of normal cone
introduced in the Appendix.}
 The map $\Be $ is uniquely characterized
by the property of being {\em asymptotic to identity} and $G\times G$-equivariant, see Corollary \ref{Bexists} (this is also an outgrowth
of an idea communicated to us by Bernstein).

The second main result of the paper addresses the question of presenting this map by an explicit correspondence. The answer is that the correspondence giving $\Be $ can be expressed in terms of the {\em inverse} to the intertwining operator (Radon transform), see Theorem \ref{AIBst}. As an application of this
result we obtain a generalization of a result of Opdam which describes the Plancherel functional restricted to a certain commutative subalgebra in the Hecke algebra.

\medskip

The paper is structured as follows. In section \ref{DCPsect} we collect necessary algebro-geometric
facts about the wonderful compactification and related varieties. Section \ref{spn}
develops the formalism of transferring functions between a $p$-adic variety and the normal bundle to its subvariety. Section \ref{alglem} contains some algebraic preliminaries. Those sections are mutually logically independent. In section \ref{consB} the results of the earlier sections are combined to get
the desired map $\Be $. In section \ref{2ndadj} second adjointness is deduced.
In section \ref{intop} we formulate some properties of the map $\Be $, including the formula
expressing it in terms of the inverse to the intertwining operator. In sections \ref{proRad} and \ref{prThm}
we prove some statements stated in section \ref{intop}. The last section \ref{Plf} contains a generalization of Opdam's formula. In the Appendix (joint with Y.~Varshavsky) we describe a version of the normal bundle construction which allows to extend some statements about De Concini - Procesi wonderful compactification to reductive groups with a nontrivial center.

We did not attempt to reach a maximal reasonable generality. In particular, we expect that
the present methods can be used to obtain a generalization of some
of our results (such as Theorem \ref{AIBst})  replacing a Borel subgroup by an arbitrary parabolic.

{\bf Acknowledgements.} We are grateful to Joseph Bernstein
 for many discussions and explanations over the years. 
We also thank  Michael Finkelberg and Yakov Varshavsky 
 for some related discussions and Vladimir Drinfeld and Jonathan Wang
 for helpful comments on the text. The authors were partly supported by the US-Israel Binational Science Foundation, NSF and ERC grants. R.B. was partly supported by the Simons Foundation
 fellowship. 

\section{Algebraic varieties related to a semi-simple group}\label{DCPsect}
 In this section
we collect algebro-geometric facts to be used in the main construction. Until the end of the section
we work over an arbitrary field 
which we denote by $F$. In this section topological notions refer to Zarisky topology.

\subsection{Standard notations}\label{basenot}
Let $G$ be  a connected split  reductive algebraic group over $F$.
 Let $\mcB$ be the space of Borel subgroups of $G$  and
 $W$ be the (abstract) Weyl group.
 Let $\mT$ be the abstract Cartan of $G$, thus $\mT$ is canonically identified with $B/U_B$
 for any $B\in \BB$, where $U_B$ is the unipotent
radical of $B$. The torus $\mT$ acts on the right on $G/U_B$ for any $B\in \BB$.

 We denote by $X_*(\mT)$ 
  the lattice of cocharacters of the torus $\mT$ and by $X^*(\mT)$ its lattice of characters.
  By $X_*(\mT)^+$, $X^*(\mT)^+$ we denote the subsets of dominant (co)characters.
 Let $\bS$ be the set of vertices of the Dynkin diagram of $G$.
  We use the bijection
between conjugacy classes of parabolic subgroups in $G$ and subsets in $\bS$ sending
a parabolic to the set of $i\in \Sigma$ such that the corresponding root subspace is in the radical of $P$.
For $I\subset \bS$ we let $P_I=L_IU_I$ denote (an arbitrary) representative of the corresponding conjugacy class.

For $i\in \bS$ we denote by $\alpha _i$ the corresponding simple root, by  $\alpha _i^\vee$ the corresponding
simple coroot and by $\bo _i$ the corresponding fundamental weight.
The abstract
Weyl group $W$ acts on $\mT$, for $i\in \bS$ denote by $s_i\in W$ the corresponding simple reflection and let $w_0$ denote the longest element of $W$.
 For a pair $B,\, B'$ of  Borel subgroups we denote by $w(B,B')\in W$ their relative position,
 with an alternative notation
$B\overset{w}{-} B'$.
We often fix a Borel subgroup $B_0\in \BB$ and denote
by $\BB_w$ the Schubert cell given by $\BB_w=\{B\in \BB\ |\ B_0
 \overset{w}{-} B\}.$

\subsection{The spaces $X$, $Y$}
Fix $B_1, B_2\in \BB$ and set $X_{B_1,B_2}=(G/U_{B_1}\times G/U_{B_2})/(B_1\cap B_2)$, where $B_1\cap B_2$ acts diagonally on the right. Given another choice
$(B_1',B_2')\in \BB^2$ in the same conjugacy class,
 an element $g\in G$ such that $B_i'=gB_ig^{-1}$ for $i=1,2$
is defined uniquely   up to
right multiplication by an element of $B_1\cap B_2$. Then the map $y\mapsto yg^{-1}$ induces isomorphisms
$G/U_{B_i}\isor G/U_{B'_i}$. The induced isomorphism $X_{B_1,B_2}\isor X_{B_1',B_2'}$ does not depend on the choice of $g$,
thus the conjugacy class of the pair $(B_1,B_2)$ defines $X_{B_1,B_2}$ uniquely up to a unique isomorphism.
Such conjugacy classes are in bijection with $W$, thus
for each $w\in W$ we get a well defined variety $X_w=X_{B_1,B_2}$,
where $w=w(B_1,B_2)$.  By construction the group $ \mT^2$ acts on varieties
$X_w$.

We will only use two extreme cases for which we fix a different notation:
$X=X_{w_0}$ and $Y=X_e$ (where $e$ is the unit element). Fixing $B\in \BB$ we get
$X=(G/U_B)^2/\mT$, $Y=(G/U_B)^2/\mT$ where the torus acts via the maps $\mT\to \mT^2$ given by
$t\mapsto (t,w_0(t))$ and $t\mapsto (t,t)$ respectively.

 We denote by $p_X:X\to \BB^2 ,\, p_Y:Y\to \BB^2$ the natural projections.

 Notice that the stabilizer in $G$ of the unit coset $(U_{B_1}, U_{B_2}) \mod (B_1\cap B_2)$
 coincides with the stabilizer of its image in $\BB^2$. It follows in particular that
 \begin{equation}\label{pXD}
 p_X^{-1}(\BB^2_0)\cong \BB_0^2\times \mT^2/\mT_\Delta \cong \BB_0^2\times \mT ,
 \end{equation}
 \begin{equation}\label{pXo}
 p_Y^{-1}(\Delta_\BB)\cong \Delta_\BB\times \mT^2/\mT_\Delta \cong \BB\times \mT,
  \end{equation}
where $\Delta_\BB\subset \BB^2$ is the diagonal,
$\BB^2_0\subset \BB^2$ is the open orbit of the diagonal $G$-action, and 
$\mT_\Delta\subset \mT^2$ is the diagonal subtorus.
The isomorphisms $\mT^2/\mT_\Delta\cong \mT$ are  given by $(t_1,t_2)\mapsto t_1^{-1} t_2$.

\begin{Ex}\label{SL2_ex}
Let $G=SL(2)$. Then $Y$ parametrizes pairs of nonzero vectors $v_1, v_2\in F^2$ modulo common dilations, while $X$ is the space of rank one $2\times 2$ matrices (notice that the symplectic form
on the two dimensional space identifies it with its dual).
 Also, $p_Y^{-1}(\Delta_\BB)$ parametrizes pairs of non-zero colinear vectors, the map $p_Y^{-1}(\Delta_\BB)\to \mT=\Gm$ sends such a pair $(v,\lambda v)$
to $\lambda$. The space $ p_X^{-1}(\BB^2_0)$ parametrizes non-nilpotent rank one matrices, the map
$ p_X^{-1}(\BB^2_0)\to \mT$ sends such a matrix to its trace.
\end{Ex}

\subsubsection{The space $X_I$}
The definition of $X$ can be generalized as follows.
Given two opposite parabolic subgroups $P$, $P^-$ with unipotent radicals $U_P$, $U_{P^-}$, an argument similar to the one above
shows that the space $(G/U_P\times G/U_{P^-})/(P\cap P^-)$ depends only on the conjugacy class of $P$.
Thus we get a variety defined uniquely up to a unique isomorphism for every
 conjugacy class of parabolic subgroups. The latter are in a bijection with
  subsets $I\subset \Sigma$ (see \ref{basenot}), for such a subset $I$ we let $X_I\cong (G/U_{P_I}\times G/U_{P_I^-})/(P_I\cap P_I^-)$  denote the corresponding variety. Thus $X=X_\Sigma$.

Similarly, set $Y_I=(G/U_{P_I}\times G/U_{P_I})/L_I$.

\begin{Rem}\label{cano}
The algebraic variety $X_I$ has been defined using a choice of the parabolic subgroup $P_I$;
however, this variety (like $G/P$ but unlike $G/U_P$ or $Y_I$) is defined canonically given the reductive group $G$.  
This is clear in view of the following: given two pairs of opposite parabolics $(P_I,P_I^-)$ and
 $(P_I', (P_I^-)')$ in the same conjugacy
class, an element conjugating  $(P_I', (P_I^-)')$ to $(P_I,P_I^-)$ is unique up to multiplication
by an element in $L$, thus the isomorphism between the corresponding quotients is canonical.
\end{Rem}

\subsubsection{Radon correspondence}
\begin{Lem}\label{fCwdef}
For any $w\in W$ there exists unique orbit $\fC_w$ of the diagonal $G^2$ action on $X\times Y$, such that

i) For any  $(x_1,x_2,y_1,y_2)\in \BB^4$  in the image of $\fC_w$ under $p_X\times p_Y$
we have $x_1\overset{w}{-} y_1$, $x_2\overset{w w_0}{-}y_2$.

ii) For any 
$x\in p_X^{-1}(\BB^2_0)$, $y\in p_Y^{-1}(\Delta_\BB)$
such that $x$, $y$ go to the unit element in $\mT$ under the projection to the second factor in the
decompositions \eqref{pXD},  \eqref{pXo}, we have $(x,y)\in \fC_w$. 
\end{Lem}

\proof It is clear that for any $x\in X$   orbits of $G^2$ on $X\times Y$ are in bijection with orbits
of the stabilizer $Stab_{G\times G}(x)$ on $Y$. If $x\in p_X^{-1}(\Delta_\BB)$, then
$Stab_{G\times G}(x)=\mT\cdot (U_{B_0}\times U_{B_0})$ for some 
 Borel subgroup 
$B_0=\mT\cdot U_{B_0}$. Thus $Stab_{G\times G}(x)$ has a unique orbit
in $p_Y^{-1}(\BB_w\times \BB_{ww_0})$ containing a point $y$ which goes 
to the unit element in $\mT$ under the projection to the second factor in \eqref{pXo}.
It is clear that this orbit is independent of the choice of $x\in p_X^{-1}(\Delta_\BB)$
provided that $x$  goes 
to the unit element in $\mT$ under the projection to the second factor in \eqref{pXD}.
\qed

We call $\fC_w$  the {\em Radon correspondence}.

\begin{Ex}\label{SL2_cor}
Let $G=SL(2)$, the spaces $X$ and $Y$ are described in Example \ref{SL2_ex}.
Then $\fC_e=\{ (m;\, v_1,v_2) \in X\times Y\ |\ m(v_2)=v_1\} $ and
$\fC_s=\{ (m;\, v_1,v_2) \in X\times Y\ |\ m^t(v_1)=v_2 \} $ where $m^t$ is the transposed
matrix with respect to the symplectic form on the two dimensional vector space and $e,\, s$
are the two elements in $W\cong \Zet /2\Zet$.
\end{Ex}

\subsubsection{Radon correspondence is closed}
Variety $X$ is well known to be quasi-affine (see e.g.
\cite[Exercise 5.5.9(2)]{Sprb}), let $\Xb_{\aff}$ be its affine closure.

\begin{Lem}\label{imageof}
The image of $\fC_w$ under the open embedding $X\times Y\to \Xb_{\aff}\times Y$ is closed.
\end{Lem}

\proof Since the subset $\fC_w\subset \Xb_{\aff}\times Y$ is $G\times G$-invariant and $G\times G$ acts
transitively on $Y$,  it is sufficient  to show that  fibers $\fC_w(y)$ of the projection $\fC_w\to Y$ are closed in
$\Xb_{\aff}$.
 The fiber over $y\in Y$  is an orbit of the stabilizer $Stab_{G\times G}(y)$. Without
 loss of generality we can assume $y\in Y$ is the image of $e\in G\times G$ and therefore
 $Stab_{G\times G}(y)=T_\Delta\cdot (U_B\times U_B)\subset B\times B$.
It is easy to see from the definition
of $\fC_w$ that the diagonal Cartan subgroup $T_\Delta$ stabilizes the point $\ti y$
in the preimage of $y$ in $\fC_w$, where
 $\ti y$ is the image of $e$. Thus this fiber is an orbit of the unipotent group
 $U_B\times U_B$.
So the Lemma follows from Kostant-Rosenlicht Theorem saying that an orbit of the action of
a unipotent algebraic group on an affine variety is closed
(see e.g. \cite[Proposition 2.4.14]{Sprb}).
\qed

\subsection{Partial compactifications of $Y$}
\label{parcomp}
Variety $Y$ is a principal bundle over $\BB^2$ with structure group $(\mT\times \mT)/\mT_\Delta\cong \mT$
where $\mT_\Delta\subset \mT\times \mT$ is the diagonal subgroup and the last isomorphism comes
from the embedding of the first factor to $\mT\times \mT$.

Recall that $X^*(\mT)^+$ is the subset of dominant weights in the weight lattice $X^*(\mT)$.
For $w\in W$ set $X^*(\mT)^+_w=w(X^*(\mT)^+)$.
We define a partial compactification $\Tb_w$ of $\mT$ by $\Tb_w:=Spec(F[X^*(\mT)^+_w])\supset
Spec(F[X^*(\mT)])=\mT$. We set $\Yb_w=Y\times ^{\mT} \Tb_w=
(Y\times \Tb_w)/\mT$, where $\mT$ acts diagonally; this is a partial compactification of $Y$.

For example, if $G=SL(2)$ then $\Yb_e$ is the total space of line bundle $\O(-1,1)$
on $\Pone \times \Pone$, while $\Yb_s$ is the total space of $\O(1,-1)$; here $e$, $s$ 
are as in Example \ref{SL2_cor}.

\subsubsection{Properness of the closure}

Let $\Cbar_w$ be the closure of $\fC_w\subset X\times Y$ in $\Xb_{\aff}\times \Yb_w$ and
 $\Cbar_w'$ be the closure of the image of $\fC_w$ in $\Xb_{\aff}\times \BB^2$
under the projection
$$Id_{\Xb_{\aff}}\times p^Y_{1,2}:\Xb_{\aff}\times Y\to \Xb_{\aff}\times \BB^2$$

\begin{Prop}\label{prprop}

a) The natural map $\Cbar_w\to \Cbar_w'$ is an isomorphism.\vspace{2mm}

b) The projection $\Cbar_w\to \Xb_{\aff}$ is proper.
\end{Prop}

\proof Since $\BB^2$ is proper, b) follows from a).

We start the proof of  a) with the following 

\begin{Lem}\label{Clnew}
a) The space $\Yb_w$ is a closed subscheme in the fiber product of total spaces of line bundles on $\BB^2$ of the form $\O(\la,-\la)$, $\la\in \Lambda^+_w$.

b)
The subscheme $\fC_w\subset X\times Y$ is the graph of sections of the pull-back of those
line bundles to a subvariety $Z\subset X\times \BB^2$; here $Z$ is the preimage under the
projection $X\times \BB^2\to \BB^4$ of the subvariety $\{(B_1,B_2,B_3,B_4)\ |\ B_1\overset{w}{-}B_3,\,
B_2\overset{w_0 w
}{-}B_4\}$.

\end{Lem}

\proof a) is clear since $\Tb_w$ is a closed subscheme in a product of affine lines,
where the torus $\Tb$ acts on each of the line by a character in $X^*(T)_w^+$. Statement (b)
amounts to the map $X\times Y\to X\times \BB^2$ inducing an isomorphism $\fC_w\to Z$.
It maps $\fC_w$ to $Z$ by Lemma \ref{fCwdef}(a); now the map $\fC_w\to Z$ is an isomorphism
since the action of $G\times G$ on both $\fC_w$ and $Z$ is easily seen to be simply transitive.
\qed

Now Proposition \ref{prprop}(a) would follow if we show that the sections 
described in Lemma \ref{Clnew}(b)
extend to the closure of $Z$
in $\Xb_{\aff}\times \BB^2$. Since $\Xb_{\aff}$ is the affine closure of a smooth quasi-affine
variety $X$,
it is normal and $\Xb_{\aff}\setminus
X \subset \Xb_{\aff}$ is
of codimension at least two, 
so any section of the line bundle on $X\times \BB^2$ extends to $\Xb_{\aff}\times \BB^2$.
Thus it suffices to show that the sections of our line bundles extend to the closure of $Z$
in $X\times \BB^2$.

Since the sections are $G\times G$ invariant and $X$ is a homogeneous space, the desired statement follows from the next

\begin{Lem}
Fix $w\in W$. For a weight $\lambda\in X^*(\mT)$ the following are equivalent

i) $\la\in X^*(\mT)^+_w$.

ii) The  line bundle $\O(\la, -\la)|_{\overline{\BB_{w_0w}}\times \overline{\BB_w}}$ has a nonzero $U\times U$
invariant section.
\end{Lem}

\proof Let us remind a description of the divisor of a $U$-invariant section of $\O(\la)$ on $\BB_w$ viewed as a rational section on $\overline{\BB_w}$, see e.g.
\cite[Proposition 1.4.5]{Brion} and references therein. Namely, for $w\in W$
components of codimension 1 in $\overline{\BB_w}\setminus \BB_w$ are in bijection with reflections $s_\alpha$, where $\alpha$ is a positive not necessarily simple  coroot, such that $\ell(ws_\alpha)=\ell(w)-1$.
The multiplicity of the corresponding component in the divisor of a $U$-invariant section of
$\O(\lambda)|_{\overline{\BB_w}}$ is then equal to $\< \lambda, \alpha \>$.
Thus those $\alpha$ for which $\ell(ws_\alpha)=\ell(w)+1$ correspond to elements of $w^{-1}(\Phi)$ which are positive roots; while those $\alpha$ for which
$\ell(w_0w s_\alpha)=\ell(w_0w)+1$ correspond to elements of $w^{-1}(\Phi)$ which are negative roots. We see that a nonzero $U$ invariant section exists both for $\O(\lambda)$ on $\overline{\BB_w}$ and for $\O(-\lambda)$ on
$\overline{\BB_{w_0w}}$ if and only if $\lambda$ is positive on $w^{-1}(\Phi)$, i.e. $\lambda\in X^*(\mT)^+_w$.

\subsection{The wonderful compactification}
\label{woncomp}
We introduce a version of De Concini -- Procesi wonderful compactification relevant for our purposes.
If $G$ is of adjoint type, we let $\Gb$ be the wonderful compactification \cite{DCP} (see also \cite{EJ},  \cite{Spr} and \cite[\S 6.1]{BrK}). If
$G$ is an arbitrary reductive group we let $G_{ad}$ denote the quotient of $G$ by its center,
$G'\subset G$ be the commutator subgroup, and consider the homomorphism $G\to (G/G')\times G_{ad}$ which is surjective and has a finite kernel. We let $\Gb$ denote  normalization of $(G/G')\times \overline{G_{ad}}$ in $G$. 

Notice that by \cite[Proposition 6.2.4]{BrK} a normal equivariant partial compactification
of $G$ with an equivariant morphism to $\Gb_{ad}$ (also known as a  {\em toroidal $G$-embedding}) is
uniquely determined by the toric variety which is the closure $\TO$ of a maximal torus $T$ of $G$ in $\Gb$.
It is easy to see (cf. \cite[Lemma 6.1.6]{BrK}) that in our case $\TO$ is the toric variety corresponding to the fan coming from
the Weyl chambers stratification.

The components of $\de G:=\Gb \setminus G$ are indexed by $\Sigma$, this easily follows e.g. from \cite[Proposition 6.2.3(ii)]{BrK}.  For $i\in \Sigma$ let $\Gb_i$ be the corresponding component

If $G$ is a product of a torus and an adjoint group,   then $\Gb$ is  smooth
  and $\de G=\Gb\setminus G$ is a divisor with normal crossings (see e.g.
  \cite[Theorem 6.1.8]{BrK}); this is not necessarily true
  in general.

 For $I\subset \Sigma$ set $\Gb_I=\bigcap \limits_{i\in I}\Gb_i$.  Let $G_I$ be the complement
 in $\Gb_I$ to the union of $\Gb_J$, $J\supsetneq I$. By convention
 $\Gb_\emptyset =\Gb$ and $G_\emptyset
 = G$.

\begin{Claim} \label{29}
a) The action of $G\times G$ on $G_I$ is transitive.

b) 
Let $P_I=L_I U_I$ be a parabolic in the conjugacy class corresponding to $I$,
 let  $P^-=L_IU_I^-$ be an opposite parabolic. Then
there exists a point in $G_I$ with stabilizer $H_I: =\{ (ul, u^- l^-) |\ u\in U_I,\,
 u^-\in U_I^-,\ l, l^-\in L_I, l^{-1}l^-\in Z(L_I)^0 \}$, where $Z(L_I)^0$ is the identity
  component in the center $Z(L_I)$ of $L_I$.
  \end{Claim}
  
  \proof a) appears, for example, as \cite[Proposition 6.2.3(iii)]{BrK}. 
  Statement b) for $G$ adjoint is \cite[Lemma 1(ii)]{Spr}, notice that in this case
  $Z(L_I)^0=Z(L_I)$.
  The general case  follows: using the adjoint case we find a point  $x_I\in \Gb$ whose stabilizer
  is a finite index subgroup in 
  $\{ (ul, u^- l^-) |\ u\in U_I,\,
 u^-\in U_I^-,\ l, l^-\in L_I, l^{-1}l^-\in Z(L_I) \}$; 
 in view of \cite[Lemma 6.1.4(i)]{BrK}  we can assume without loss of generality that $x_I$ lies in the closure $\TO$
 of a maximal torus $T\subset L_I$, then the above description of $\TO$ shows that the stabilizer of  $x_I$ in $T$ is connected, thus
 that the stabilizer of $x_I$ in $G\times G$ is connected and coincides with $H_I$. 
  
  See also   
\cite[Theorem 2.22, Proposition 2.25]{EJ} for the case of an adjoint group over a field of characteristic zero. 
    \qed

\subsubsection{Normal bundles to strata}\label{normsec}
The isomorphism $G_I=(G\times G)/H_I$ shows that we have a canonical map $X_I\to G_I$ which is
a principal bundle with the structure group $Z(L_I)^0$.

For a locally closed smooth subvariety $Z$ in a smooth variety $M$ let $N_M(Z)$ denote
the normal bundle. More generally, in the situation described in the Appendix, section \ref{qnc},
we let $N_M(Z)$ denote the quasi-normal cone in the sense of \ref{qnc}.

If $G$ is adjoint,
the normal bundle $N_{\Gb}(\Gb_I)$ splits canonically as a sum of line bundles $N_{\Gb}(\Gb_i)|_{\Gb_I}$, $i\in I$. We have an action of $Z(L_I)$ on  $N_{\Gb}(\Gb_I)$ such that
the action on $N_{\Gb}(\Gb_i)|_{\Gb_I}$, $i\in I$ is by the character $\alpha_i$. We will need 
a generalization of this observation to an arbitrary $G$.

Let 
$A_I$ denote the closure of $Z(L_I)^0$ in the $\Tb_{\aff}:=Spec(F[\Lambda^+])$.
 This is a toric variety for the torus $Z(L_I)^0$.

\begin{Claim}\label{normbdl}
a) The fibration $N_{\Gb}(G_I)\to G_I$ is canonically isomorphic to the bundle with fiber
$A_I$ 
 associated to the principal $Z(L_I)^0$ bundle $X_I\to G_I$.


b) There exists a canonical $G\times G\times Z(L_I)$ equivariant open embedding $X_I\to N_{\Gb}(G_I)$.
Its image is the complement to the divisor $\bigcup\limits_{i\in  I} N_{\Gb_i}(G_I)$.
\end{Claim}

\proof 
(b) follows from (a), while (a) follows from Theorem \ref{qnrml}(c) in view of Example
\ref{exl1}. Alternatively,
part (a) can be deduced
from \cite[Proposition 6.2.3(i)]{BrK}  which shows that 
the variety $\Gb$ is locally isomorphic to the product
of the affine toric variety $\Tb_{\aff}$ with a smooth variety, the stratification
of $\Gb$ by $G\times G$ orbits $G_I$ corresponds to the stratification by $\mT$ orbits
on $\Tb_{\aff}$. \qed


\begin{Rem}

The following description of $X_I$, though not used explicitly in this paper, is closely related
to Claim \ref{normbdl}. The space $X_I$ is quasi-affine, it is the dense
 $G\times G$ orbit in its affine closure $\Xb_{I,\aff}$. Thus $X_I$ can be reconstructed from the algebra
 of global regular functions $\OO_{gl}(X_I)$. Assume for simplicity that $F$ is of characteristic zero.
 Then the space $\OO_{gl}(G)$
 is isomorphic as a $G\times G$-module to: $\oplusl_{\lambda\in X^*(\mT)^+} E_\lambda$,
 where
$E_\lambda= V_\lambda \otimes V_\lambda^*$ for the representation $V_\lambda$ with highest weight $\lambda$. Let $m_{\lambda,\mu}^\nu:E_\lambda\otimes E_\mu\to E_\nu$ be the corresponding component of multiplication in $\OO_{gl}(G)$. Then $\OO_{gl}(X_I)=\OO_{gl}(G)$ as $G\times G$ modules and
multiplication map in $\OO_{gl}(X_I)$ is the sum of maps 
$m_{\lambda,\mu}^\nu(I):E_\lambda\otimes E_\mu\to E_\nu$, where
 $m_{\lambda,\mu}^\nu (I)=m_{\lambda,\mu}^\nu$
when $\lambda+\mu-\nu$ is trivial on $Z(L_I)$ and $m_{\lambda,\mu}^\nu (I)=0$ otherwise.
\end{Rem}

\begin{Rem}
The correspondence $\fC_w$ can also be described in terms of geometry of the wonderful compactification $\Gb$. Namely, let $\Gamma \subset G\times \BB^2$ be the graph of the action of $G$ on $\BB$, i.e. it is given by: $\Gamma=\{(g,B_1,B_2)\ g(B_1)=B_2\}$.
Let $\bar{\Gamma}$ be the closure of $\Gamma$ in $\Gb\times \BB^2$. Let
$Z\cong \BB^2 \subset \Gb$ be the closed $G^2$ orbit.
According to \cite{Br},
\cite{BrP}, the subset
$$(Z\times \BB^2)_w:=\{(B_1, B_2; B_3, B_4)\ |\ B_1\overset{w}{-}B_3,
B_2\overset{w_0w}{-}B_4\}$$
is an open smooth subscheme in $\bar{\Gamma}\cap Z\times \BB^2$.
We claim that $\fC_w$ is canonically isomorphic to an open part in the quasi-normal
cone (see Appendix for a definition of this term) $N_{\bar{\Gamma}}((Z\times \BB^2)_w)$. The projection $\fC_w\to X$ is the restriction
of the natural map (the normal differential of the projection, i.e. the composition of the differential 
with projection to the normal bundle) $N_{\bar{\Gamma}}((Z\times \BB^2)_w)
\to N_{\Gb}(Z)\supset X$. The projection $\fC_w\to Y$ can be described as follows.
Recall that there exists a canonical action map $a:\Gamma=G\times \BB\to Y$, $a:(g,B)\to (\tilde{B}, g(\tilde{B})\mod T)$, where $\tilde{B}$ is an arbitrary
lifting of $B\in \BB$ to a point in $G/U$. Let us view $a$ as a rational map
$\bar{\Gamma}\to \Yb_w$. We claim that $(Z\times \BB^2)_w$ is contained in the domain
of definition of this map, and $a:(Z\times \BB^2)_w\to \BB^2\subset \Yb_w$. Thus
the normal  differential of $a$ gives a map $N_{\bar{\Gamma}}((Z\times \BB^2)_w)\to
N_{\BB^2}(\Yb_w)=\Yb_w$, which restricts to the projection $\fC_w\to Y$.

We neither prove nor use these statements.
\end{Rem}

\subsubsection{Quotient by the $U_P$ action}
Let $P=P_I=L_I U_I$ be a parabolic subgroup whose conjugacy class corresponds to $I\subset \Sigma$, let $Z_I=Z(L_I)^0$ denote the neutral component in the center of $L_I$.

The map $G/U_I\to G/(U_I Z(L_I)^0)$ is a principal $Z(L_I)^0$ bundle.
We let $\overline{G/U_I}:=(A_I\times G/U_I)/Z_I$ 
be the associated bundle over $G/(U_I Z_I)$ with fiber $A_I$.

In view of Claim \ref{29}(b) we have a canonical projection $G_I\to G/P_I^-$;
let $G_I^0\subset G_I$ be the preimage of the open $U_I$ orbit under that projection.

Consider the closure of the subgroup $Z_I$ in $\Gb$. This is a toric variety for the torus $Z_I$.

\begin{Claim}\label{qUP}  
a) The closure of $Z_I$ intersects $G_I^0$ at a unique point, which we denote by $z$. 

b) Let $\tilde Z_I$ be the open subvariety in this toric variety which is the union of all $Z_I$-orbits
whose closure contains $z$. Let $\Gb^0(I)$ denote the image of $\tilde Z_I$ under the action of
$G\times U_I$.
Then the subset $\Gb^0(I)\subset \Gb$ is an open subvariety, such that

i) $\Gb^0(I)$ is invariant under the left action of $G$ and the right action of $U_I$, the action of $U_I$ on $\Gb^0(I)$
is free.

ii) the quotient $\Gb^0(I)/U_I$ is isomorphic to $\overline{G/U_I}$.

iii) We have $G\subset \Gb^0(I)$, $G_I^0\subset \Gb^0(I)$ and the induced embeddings of the quotient
spaces coincide respectively with the tautological embedding $G/U_I\to \overline{G/U_I}$
and the embedding  $G/(U_I Z_I)\to \overline{G/U_I}$ induced by the embedding $\{0\}\to
A_I$.

\end{Claim}

\proof Fix a maximal torus $T\supset Z_I$ and Borel subgroups $B=TU\subset P_I$ and $B^-=TU^-\subset
P_I^-$. 
According to \cite[Proposition 6.2.3(i)]{BrK} there exists an open subvariety $\fU\subset \Gb$
 (denoted by $X_0$ in {\em loc. cit.}) such that the action of $U^-\times U$ on $\fU$ is free, where $U^-$ acts on the left 
 and $U$ acts on the right. Furthermore, $\fU\cong U^-\times \Tb_{\aff}
 \times U$  (where $\Tb_{\aff}=A_\Sigma$ is as in the proof of Claim \ref{normbdl}; in {\em loc. cit.} it is denoted by $X_0'$).
By  \cite[Proposition 6.2.3(ii)]{BrK} the intersection  $\Tb_{\aff}
\cap \Gb_I$ is a $T$ orbit,
to be denoted by $\Tb_{\aff} (I)$. The closure of $Z_I$ in $\Tb_{\aff}$
(recall that it is denoted by $A_I$) intersects the $T$-orbit $\Tb_{\aff}(J)$ iff $J\supseteq I$. Also
the intersection of this closure with $\Tb_{\aff}(I)$ consists of a single point $z$. 
This point $z\in G_I$ lies on a free $U^-\times U$ orbit, hence it lies in $\Gb^0(I)$. Conversely,
any other point $ z'$ in the intersection of $G_I$ with the closure of $Z_I$ 
is stabilized by the diagonal action of $L_I$, as well as by a subgroup which is conjugate to but different from $U_I\times U_I^-$ and is normalized by $L_I$, hence $z'$ does not lie in $G_I^0$.  This proves (a). 
Notice also that we have shown that $A_I\cong \tilde Z_I$. 


To check that $\Gb_0(I)$ is open it is enough to show that it contains a neighborhood of $Z_I$.
For such a neighborhood we can take the image of $\tilde Z_I$ under the action of $U^-\times T\times U$: on the one hand it coincides with the 
$\cupl _{J\supseteq I} U^-\Tb_{\aff}(J) U$ which is open in $\fU$, hence in $\Gb$; on the other hand it is contained in $\Gb^0(I)$ since $B\subset L_I\cdot U_I$ and $L_I$ commutes with $Z_I$.

Since $\tilde Z_I$ is contained in  $\fU$, it lies in the union of free (right) $U_I$ orbits, so the 
right action of $U_I$ on $\Gb_0(I)$ is free. The image of the map $A_I\cong \tilde Z_I \mapsto 
\Gb_0(I)/U_I$ is clearly invariant under $U_I Z_I$ and the action of $U_I$ on this image is trivial,
while the action of $Z_I$ is the canonical one. Thus we get a surjective map 
$\overline{G/U_I}\to \Gb^0(I)/U_I$. The stabilizer of any point in $\tilde Z_I$ under the left action of 
$G$ is contained in $U_I$, hence the stabilizer of any point in the image of $\tilde Z_I$
under the map to $ \Gb^0(I)/U_I$ equals $U_I$, this implies that the map $\overline{G/U_I}\to \Gb^0(I)/U_I$ is an isomorphism. Property (iii) is clear from the construction. \qed

\begin{Ex}
Let $G=PGL(2)$, thus $\Gb=\bbP(End(V))$ where $V=k^2$ is the two dimensional vector space. Let $P_I=B_I$
be the group of upper triangular matrices. Then $\Gb^0(I)$ is the projectivization of the set of matrices with a nonzero second column. The quotient $\Gb^0(I)/U_I$
 maps to $\bbP^1$ by the map sending a matrix
to the line of its second column, this identifies $\Gb^0(I)/U_I$ with the total space of the line bundle
$\OO_{\bbP^1}(2)$.
\end{Ex}

\subsubsection{Compactified Bruhat cells} \label{Brcl}
Fix $B\in \BB$. For $w\in W$  let $G_w\subset G$ denote the corresponding $B\times B$ orbit.
 Let $\overline{G_w}$  denote the closure of $G_w$ in $\Gb$.

To unburden notations we state the next Claim for a semi-simple group $G$, the answer in the general case differs by replacing $\BB^2$ and its subvarieties by their product with the torus 
$G/G'$.

 \begin{Claim}\label{BrclClaim}
a)  The components of the intersection of $\overline{G_w}$ with the closed stratum $\BB^2
\subset \Gb$ are in bijection
with pairs $w_1$, $w_2$
such that
 $\ell (w)+\ell(w_0)=\ell(w_1)+\ell(w_2)$.

b)
Given such a pair
$(w_1,w_2)$, the smooth part of the corresponding component contains the $B\times B$ orbits
$\BB_{w_1}\times \BB_{w_2}$.

c) The
open part of the quasi-normal cone to the smooth locus of the corresponding component
is identified with $X_{w_1,w_2}$ where $X_{w_1,w_2}$ is the corresponding
 $B\times B$ orbit in $X$.
\end{Claim}

\proof If $G$ is adjoint then (a,b) follow from \cite[Theorem 2.1(ii)]{Br}. In view
of Claim \ref{normbdl}, part (c) in this case follows from \cite[Theorem 2.1(i)]{Br}
 which implies that the normal bundle to $\overline{G_w}\cap \BB^2$ in
 $\overline{G_w}$ is isomorphic to the restriction of the normal bundle to $\BB^2$ in $\Gb$.
 The general case follows since $\overline{G_w}$ is easily shown to coincide with the preimage of
 of $\overline{(G_{ad})_w}$, this  implies in particular that the quasi-normal cone to 
 the smooth part of $\overline{G_w}\cap \BB^2$ is the fiber product of the quasi-normal cone to
 $\BB^2$ in $\Gb$ by the smooth part of $\overline{G_w}\cap \BB^2$. \qed
 
\begin{Cor}\label{BrclCor} For $w$, $w_1$, $w_2$ as in Claim \ref{BrclClaim} 
we have a canonical isomorphism of $\mT$ torsors: $U_B\bs G_w/U_B \cong U_B\bs
X_{w_1,w_2}/B_U$. 
\end{Cor}

\proof It easy to see that for a smooth toric variety $C$ the normal bundle
to a torus orbit is canonically identified with an open subvariety in $C$.
The same is true for a not necessarily smooth normal toric variety with a fixed finite 
equivariant map to a smooth one, where instead of the normal bundle we consider
the quasi-normal cone introduced in the Appendix. 

One can apply this
to the closure $C_w$ of a generic orbit of a maximal torus $T$ acting on $G_w$ by left translations
and the $T$ orbit $C_{w_1,w_2}:=C\cap \BB_{w_1}\times \BB_{w_2}$. The open orbit of $T$
on $C_w$ maps isomorphically to $U_B\bs G_w/U_B$, while the open 
orbit in $N_{C_{w_1,w_2}}(C)$ maps isomorphically to $U_B\bs X_{w_1,w_2}/U_B$.
Thus the sought for isomorphism is obtained by
restricting the isomorphism described in the previous paragraph to the open
$T$ orbit. 
\qed

\subsection{More on partial compactifications of $X$}

\subsubsection{Torus closure}

Let $\Tb:=\Tb_1$ denote  the partial compactification of the abstract Cartan $\mT$ attached
to the unit element in $W $ by the construction of section  \ref{parcomp}.

\begin{Claim}\label{clXaff}
a) The closure in $\Xb_{\aff}$ of any orbit of the abstract Cartan $\mT$ acting on $X$ is isomorphic to
$\Tb$.

b) For any $w_1,w_2$ there exists
a $\mT$ equivariant isomorphism $\mT\cong X_{w_1,w_2}/U_B^2$.
 The corresponding projection $X_{w_1,w_2}\to \mT=X_{w_1,w_2}/U_B^2$ extends to a regular map
$\overline{X_{w_1,w_2}}\to \Tb$ where $\overline{X_{w_1,w_2}}$ denotes the closure in $\Xb_{\aff}$.

\end{Claim}

\proof It is clear that the weights of $\mT$ acting on the space of regular functions on $X$
are precisely the dominant weights; also, given such a weight $\la$ and an orbit of 
$\mT$ there exists a function transforming by the character $\la$ with a nonzero restriction to the
orbit. This implies part (a). The first statement in part (b) is clear, the second one amounts to existence for every dominant weight $\la$ of a section of the line bundle $\O(\la,\la)$ on $\BB^2$ whose restriction to 
 $\BB_{w_1}\times \BB_{w_2}$ is nonzero and $U_B^2$ invariant. Since $U_B^2$ is unipotent
 it suffices to see existence of a section with a nonzero restriction which is clear. \qed

\subsubsection{Embedding into the matrix space and its compactification}
\label{embinos}

For a vector space $\VV$ let $\bbP(\VV)$ denote the projectivization of $\VV$, let $\bbPt(\VV)=\bbP(\VV\oplus F)=
\bbP(\VV)\cup \VV$ be its projective compactification, and set $\VV^0=\VV\setminus \{0\}$,
$\bbPt(\VV)^0=\bbPt(\VV)\setminus \{0\}$.

We set $\Xb=N_{\Gb}(\Gb_\Sigma)$; in view of Claim \ref{normbdl}
it is identified with the fibration over $G_\Sigma$ associated to the $\mT$-bundle
$X \to  G_\Sigma$ and the $\mT$ space $\Tb_{\aff}$. In other words, $\Xb$ is the relative
spectrum of the sheaf $\bigoplus\limits_{\la\in \La^+} \O_\BB(\la)\boxtimes
\O_\BB(\check{\la})$ of commutative rings on $\BB^2$ where  $\check{\la}=-w_0(\la)$ is the dual 
dominant weight.

\begin{Claim}\label{clinEnd}

a) Let $M$ be a $G$-module. Then there exists a canonical map $\rho_M:X\to End(M)$.
This map extends to a map $\overline{\rho_M}:\Xb\to \bbPt(End(M))^0$.

b) There exists a finite collection of $G$-modules $M_i$, such that the map
$\prod \overline{\rho_{M_i}}:\Xb\to \prod \bbPt(End(M_i))^0$ is a closed embedding.
This embedding sends the zero section of $N_{\Gb}(\Gb_\Sigma)=\Xb$ to $\prod \bbP(End(M_i))$.
\end{Claim}

\proof Fix a pair of opposite Borels $B=TU$, $B^-=TU^-$, so that $X=(G/U\times G/U^-)/T$.
Let $\la_i$ be the set of highest weights of $M$ (i.e. maximal elements in the set of weights of $M$)
and let $E_{\la_i}\in End(M)$ be  $T$-invariant projection to the corresponding weight space.
The operator $E=\sum \limits_i E_{\la_i}$ is invariant under the left action of $U$, the right action
of $U^-$ and the diagonal action of $T$, thus we get a map $\rho_M:X\to End(M)$ sending the unit coset
to $E$. Since any two pairs of opposite Borel subgroups are conjugate by an element of $G$ which is defined uniquely up to multiplication by an element in the Cartan subgroup, the map $X\to 
End(M)$ does not depend on the auxiliary choice.
The extension to map $\overline{\rho_M}: \Xb\to \bbPt(End(M))$ comes from the standard highest weight morphism of $G^2$ equivariant bundles over $\BB^2$: $M^*\otimes M \otimes \O_{\BB\times \BB}
\to \bigoplus\limits_i \O(\la_i, \check{\la_i})$; this proves part (a).
According to \cite[Lemma 6.1.1]{BrK} there exists a module $M_0$ for which the map
$\BB\times \BB\to \bbP(End(M_0))$ induced by $\rho_{M_0}$ is a closed embedding.
We can choose modules $M_1,\dots, M_k$ and a highest weight $\la_i$ of $M_i$ so that
the weights $\la_1,\dots, \la_k$ generate the lattice of dominant weights. 
Then the map $\prod\limits_{i=0}^k \overline{\rho_{M_i}}$ is easily seen to be a closed
embedding. The last property stated in (b) is clear from the construction. \qed

\section{Co-specializaton from  normal bundle for functions on a $p$-adic manifold}\label{spn}
From now on $F$ is a non-Archimedian local field with the ring of integers $O\subset F$ and a uniformizer $\pi$. From now on topological notions are in reference to  the $F$-topology on the space of $F$-points of an algebraic variety over $F$ or its subspaces.

For a totally disconnected space $X$ we let $\SS(X)$ denote the space of  locally constant compactly supported functions on $X$. 

Let $\WW$ be a smooth analytic variety over a non-archimedian local field $F, D\subset \WW$
an open subset such that the complement $\WW\setminus D$ is a union $S=\cup _{i\in \bS}S_i$ of smooth divisors
$ S_i$ with normal crossing.
 For $I\subset \bS$ we define  $S_I=\cap_{i\in I} S_i$,  $S_I^0:=S_I-\cup _{j\in \bS -I}S_{I\cup j}$ and denote 
 by $r_I:N _I\to S_I$ the normal bundle to $S_I$.
For any $J\subset I$ we denote by $N_I^J\subset N _I$ the normal bundle to $S_I$ in $S_J$.
 Locally in $S_I$ we can identify  $N_I$ with the product $\mA ^{ I}\times S_I$. 


The following is immediate.
\begin{Claim}
\label{normal} a) We have a canonical direct sum decomposition
$$N_I=\oplus _{i\in I} N_I^{I-\{ i\}}$$

b)  $r_J^* (N_I^J)$
is the  normal bundle to $r_J^{-1}(S_I)$ in $N_J$. 

c) The complement $N_I-r_I^{-1}(S^0_I)$ is a union of smooth divisors 
$r_I^{-1}(S^{I\cup j}),j\in \bS \setminus I$ with
normal crossing.
\end{Claim}
\begin{Def}
  Let $U\subset \WW$ be an open subset. We say that an  analytic  open embedding
${\tau _I} :U\to N_I$  is {\it admissible} if \vspace{2mm}

${\tau _I} _{|S_I\cap U}=Id$  \vspace{2mm}

${\tau _I} (U\cap S_J) \subset N_I^J$ for $J\subset I$ and \vspace{2mm}

The normal component of $d{\tau _I} _{|S_I\cap U}$ 
equals the canonical projection $T\WW|_{S_I}\to N_I$.
\vspace{2mm}

\end{Def}
 Any admissible map ${\tau _I} :U\to N_I$ defines  an embedding ${\tau _I} ^* _V:\SS(V)\to \SS (U)$
for any   open subset $V$ of ${\tau _I} (U)\subset N_I$.

We denote by $(a,y)\to ay $ the natural action of the group $\mG _m^I$ on the
bundle $\un p_I: N_I\to  S_I$ and define $X_I\subset N_I$ as the open
subvariety of points with trivial stabilizer; thus $X_I=N_I\setminus \cup_{j\in I} N_I^{\{ j\}}$.

\begin{Def}\label{T} a) For any  $\bl \in Hom (\mG _m , \mG _m^I)=\mZ ^I$ we write
$z_\bl :=\bl (\pi^{-1}) \in (F^\times)^I$, $z_I:=z_{(1,\dots,1)}$.  \vspace{2mm}

b) For any $\bl \in \mZ ^I$
we define an operator $T^\bl \in End (\SS(N _I))$ by
$$T^\bl (f)(y) :=f(z_\bl y)$$

\end{Def}

For $\la=(\la_i)$, $\mu=(\mu_i)\in \Zet^I$ 
we will say that $\bl \geq \mu$ if $ \bl _i\geq \mu _i$ for all
$i\in I$.

\begin{Lem}\label{well_defined}  For any $f\in  \SS(N _I)$ and a pair
${\tau _I}  _1,\, {\tau _I} _2:U\to N_I$  of  admissible maps
there exists $\bl _0\in \mZ^I$ such that $supp (T^\mu (f))\subset {\tau _I}  _1(U)\cap {\tau _I}  _2(U)$ and
 ${\tau _I} _1^\star  (T^\mu (f))={\tau _I} _2^\star  (T^\mu (f))$ for all $\mu \geq \bl_0$.

\end{Lem}

{\em Proof.} The following obvious Claim implies the lemma.\vspace{2mm}

Let $S$ be an analytic $F$-variety, $\ti S:=F^n\times S$ and let
${\tau} : \ti S\to \ti S$ be an analytic morphism such that: \vspace{2mm}


a) ${\tau}| _{\{ 0\} \times S}=Id$.\vspace{2mm}

b) $d{\tau }| _{\{ 0\} \times s}$ induces the identity map $F^n\to F^n$ for every $s\in S$.\vspace{2mm}

Let $t\in (F^\times)^n$ be such that $|t_i|<1$.
Consider the sequence of morphsims
$$\phi _n: \ti S \to \ti S, \ \phi _n (\ti s):=t^{-n}{\tau _I} (t^n\ti s), \ {\mathrm{where}}\ \  t(v,s):=(tv,s).$$

\begin{Claim}\label{conv} We have $\phi_n(x)\to x$ as $n\to \infty$
for every $x\in \ti S$.  Moreover $\phi_n\to Id$ uniformly on every compact subset
in $\ti S$.
\end{Claim}

\begin{Def}
For a set $\fS$ consider the set of 
partially defined maps $\Zet^I\to \fS$ whose domain contains a set of the form $\{\la \ | \ \la\geq \la_0\}$ for some $\la_0\in \Zet^I$. We say that two such maps
agree asymptotically if they
coincide on $\{\la \ | \ \la\geq \la_0\}$ for  $\la_0\in \Zet^I$.
This defines an equivalence relation on the set of such maps, let $\fS_{\as}(I)$ denote
the set of equivalence classes.
\end{Def}

We restate the result of  Lemma  \ref{well_defined}
by saying that  we have a well-defined {\it asymptotic embedding}
of $$\SS(N _I)\to \SS(\WW)_{\as}(I)
,\ \ \ \ f\mapsto (\la\mapsto  {\tau _I} ^\star  (T^\bl (f))).$$

\medskip

Given two subsets $J\subset I\subset \bS$ we can apply Lemma
 \ref{well_defined} to the case considered in Claim \ref{normal}
to obtain the  asymptotic embedding
$$j_I^J : \SS(N _I)\to \SS(N_J)_{\as}(I\setminus J),\ \ \  f\mapsto  (\la\mapsto {\tau _I} ^\star  (T^{ \ti \bl} (f)).$$
The following is immediate.

\begin{Claim}\label{compatible}
For any $f\in  \SS(N _I)$ we have
$$j_J  ( j_I^J(f)( \ti \bl))(\mu)=j_I (f) (\ti \bl +\mu )$$
for $\ti \bl \in  \mZ ^I/\mZ ^J,\mu \in \mZ ^J,\ti \bl ,\mu \gg 0$ where we use the canonical isomorphism
$ \mZ ^I=\mZ ^I/\mZ ^J\oplus \mZ ^J$
\end{Claim}

The following two statements follow easily from the definitions.\vspace{2mm}

Suppose that  an analytic unimodular $F$-group $H$   acts freely on $\WW$ preserving irreducible components of $S=\WW\setminus D$.
Then for any strata $S^0_I$ we have a natural action of $H$ on $N_I$ and on $N^0_I$. We denote by
$\pi :\WW\to \bar \WW :=\WW/H$  the natural projection. Then $\pi $ defines projection $\pi _I:N_I\to \bar N _I$
where $\bar N _I$ is the normal bundle to $ \bar S_I^0:=S_I^0/H$ in $\bar \WW$. We define push-forwards
$\pi _\star :\SS (\WW)\to \SS (\bar \WW),{\pi _I}_\star :\SS (N _I)\to \SS (\bar N _I )$ as integration along $H$.

\begin{Claim}\label{compfac} The asymptotic embeddings
$$j_I:\SS(N_I)\to \SS (U)_{\as}(I),\ \ \ \ \ \  \bar  j_I :\SS(\bar N_I)\to \SS (\bar U)_{\as}(I)$$ are  compatible with
push-forwards $\pi _\star ,{\pi _I}_\star$.
\end{Claim}

\begin{Claim}\label{loccl}
Let $Z\subset \WW$ be a locally closed smooth subvariety which intersects each of $S_I$
transversely. Then the asymptotic embeddings for  $\WW$ and  $Z$ are compatible under the restriction map.
\end{Claim}

 Let  $U_I\subset \WW,I\subset \bS$  be open neighborhoods of $S_I$  and
$$\tau _I:U_I\to N_I,\, I\subset \bS ,\ \ \ \ \tau _I ^J :\ti N^I_J\to  r_J^\star (N_J^I),\, J\subset I$$
be a family of  admissible  maps (cf. Claim \ref{normal} b) where $\ti N^I_J\subset N_J$ are open neighborhoods of
$r_J^{-1}(S_I)$.

\begin{Def}\label{family}We say that a  family $(U_I,\tau _I,\tau _I ^J)$ is {\it admissible} if for any pair
$\emp \subset J\subset I\subset  \bS$  there exists an open  neighborhood $U'\subset U$ of $S_I$ in $\WW$  such that
the restriction of $\tau _I ^J \circ \tau _J$ on $U'$ coincides with $\tau _I$.

\end{Def}
The following result is easy to prove by induction on $|\bS |$.

\begin{Claim}\label{adm} An admissible family exists for any $\WW$, $D$ as above.
\end{Claim}

\subsection{The almost homogeneous case} In this section we suppose that the variety 
$\WW$ (hence also $S=\WW\setminus D$) is compact, that
 an analytic unimodular $F$-group $H$  acts on $\WW$ preserving irreducible components of $S$ and that
 the action on $S_I^0$ is transitive for all $I\subset \bS$. Let $K\subset H$ be an open compact subgroup
and  $(U_I,\tau _I,\tau _I ^J)$ be an admissible  family as in Definition \ref{family}.  \vspace{2mm}


\begin{Lem}\label{limit} For any $s\in S_I$ and a compact subset $C$ of $H$ there exists an open neighborhood $U$ of
$s$ in $\WW$ such that $\tau _I(cKu)=cK\tau _I(u)$ for all $c\in C,u\in U$.
\end{Lem}

\proof  We can assume without loss of generality that $U=cK$ for some $c\in H$.
Then the statement follows by applying Lemma \ref{well_defined} to
 ${\tau_I}_1=\tau_I$,  ${\tau_I}_2
=\tau_I\circ c$ and $f=\delta_{Ku}$ for some $u\in U$. 
\qed

\begin{Prop}\label{compatibility} For any compact subset $C$ of $H$ there exists an open neighborhood $V$
of $S$ in $\WW$, such that $V\subset \cup U_I$ and $\tau _I (h K v)=hK \tau _I(v)$ for all $h\in C$,
$v\in V$,  $I\subset \bS$.
\end{Prop}
{\bf Proof}. By induction in $n$ we construct an open set $V_n$ such that 
the required equality holds for all $I$ with $|\bS -I| \leq n$.
\vspace{2mm}

Consider first the case $I=\bS$. Then the variety $S^0_I$ is compact and therefore  is a finite union of $K$-orbits.
Therefore there exist a $K$-invariant open neighborhood $U$ of $S_I$ in $N_I$ and a finite number of points
$u_a\in U,a\in A$ such that  any $u\in U\cap D$ can be written in the form
$$u=kT_\bl u_a,\ \ \ \ a\in A,\ k\in K,\ \bl \in {\mZ ^I}_+.$$
The statement now follows  from Lemma \ref{well_defined}.\vspace{2mm}

Assume now that the statement is known for all $J\subset \bS$ containing $I$. For any $j\in \bS -I$ we choose an
open neighborhood $\ti N^I_{I\cup j}$ of $r_{I\cup j}^{-1}(S_I)$ satisfying the 
requirement of Definition  \ref{family}. Then there exist a
$K$-invariant open neighborhood $U$ of $S_I$ in $N_I$ and a finite number of points
$u_a\in U,a\in A$ such that  any $u\in U\cap D-\cup _{j\in \bS -I}\ti N^I_{I\cup j}$ can be written in the form
$$u=kT_\bl u_a,\ \ \ \ \ a\in A,\  k\in K,\ \bl \in {\mZ ^I}_+.$$
The statement  follows now from Lemma \ref{well_defined} and the definition of the admissible family.
\qed

\subsection{The singular case}\label{sc} We now drop the assumption
 that the analytic variety $\WW$ is
smooth and the divisor $S\subset \WW$ is a divisor with normal crossing. Instead we assume that
each stratum $S_I$ is {\em well approximated} in the sense of the Appendix. The definition
of an admissible map applicable to this case is given in the Appendix.
Other definitions
and results of the present section carry over to this case {\em mutatis mutandis}.

\section{An algebraic Lemma}\label{alglem}

Let $A$ be a 
 Noetherian $\mC$-algebra, $M$ a finitely generated $A$-module, $T$ an
automorphism of $M$ as an $A$-module.







\begin{Def} \label{compadef}
Let $M$, $T$ be as above and  $M^0\subset M$ be a $T^{-1}$-invariant $\mC$-subspace  such that
$M=\cup _{r>0}T^rM^0$. Let $N$ be an   $A$-module.
 A $\mC$-linear map $a:M^0\to N$ will be called {\it $A$-compatible} if for any
$m\in M$,  $h\in A$ there exists an integer
$N(m,h)>0$ such that $ha(T^{-n}m)=a(hT^{-n}m)$ for $n>N(m,h)$.
\end{Def}

\begin{Lem}\label{ext}
For  any $A$-compatible $\mC$-linear map $a:M^0\to N$ there exists
 unique  $A$-morphism
$\ti a:M\to N$ such that 
for any $m\in M$ we have  $\ti a (T^{-n}m)=a(T^{-n}m)$ for $n\gg 0$. 
\end{Lem}

\proof  Fix a finite set
   of generators $m_i$
  of $M$ over $A$.
For some $n>0$ we have $T^{-n}m_i\in M^0$.
If $\ti a$ satisfies the assumption of the Lemma, then for large $n$ we have:
 $\ti a(T^{-n}m_i)
=a(T^{-n}m_i)$; however, for any $n$ the elements $T^{-n}m_i$ generate $M$, so uniqueness
is clear. 


We now proceed to show existence. Assume without loss of generality that the generators
$m_i$ lie in $M_0$. 
Let $\ti M\cong A^s$ be the free $A$ module
 with free generators denoted $\ti m_1, \dots \ti m_s$,
and $\Pi:\ti M\to M$ be the surjection sending $\ti m_i$ to $m_i$. 
Let $\tau:\ti M\to \ti M$ be an endomorphism lifting 
$T^{-1}$.

For  $n\geq 0$  define a homomorphism $\alpha_n:
\ti M\to N$ by requiring  that $\alpha_n(\ti m_i)=a (T^{-n}m_i)$.
Given $x\in \ti M$ the $A$-compatibility condition implies  that $$\alpha _n(x) = a(T^{-n}\Pi(x))$$
for large $n$.

Applying this to $x_i=\tau(\ti m_i)$  we conclude that
 there exists $n_0>0$ such that
$\alpha_{n+1}=\alpha_n \circ \tau$ when $n\geq n_0$.
Set $\alpha=\alpha_{n_0}$.

Given $x\in \ti M$ we have 
\begin{equation}\label{aln}
\alpha(\tau^n x) = a(T^{-n-n_0} \Pi(x))
\end{equation}
for large $n$.

 Since $\tau$ induces an invertible map on $\ti M/Ker(\Pi)=M$, it follows
that for $n>0$ we have: $\Pi(x)=0 \iff \Pi \tau^n(x)=0$, hence
  $\tau^n:Ker(\Pi)\twoheadrightarrow (Ker (\Pi )\cap Im(\tau^n))$.

Since $A$ is Noetherian, $Ker(\Pi)$  is finitely generated. 
Fix a finite set of generators,
and pick $n$ for which \eqref{aln} holds when $x$ belongs to this set of generators.
Then $\alpha$ vanishes on $Ker(\Pi)\cap Im(\tau^n)=\tau^n(Ker(\Pi))$. Since $\Pi|_{Im(\tau^n)}$
 is onto,
the homomorphism $\alpha|_{Im(\tau^n)}$ factors through a homomorphism
$\ti a':M\to N$.  Now set $\ti a = \ti a'\circ T^{n_0}$, then \eqref{aln} implies that 
$\ti a$ satisfies the required condition.

\section{Construction and properties of the maps $\Be _I$}\label{consB} 
We modify the meaning of notation: from now on $G,\, B,\, \BB,\, X$ etc. will denote the groups/spaces of $F$-points
of algebraic groups or varieties considered in section \ref{DCPsect}; these are equipped with the $p$-adic topology.

We apply the construction of section \ref{spn} to $\WW=\Gb$ being the wonderful compactification, with $D=G$.
If $G$ is adjoint this is a particular case of the situation introduced at
 the beginning
 of section \ref{spn}, otherwise,
one should apply considerations of subsection \ref{sc}.
 We fix an admissible system $\tau_I$.

\subsection{Definition of $\Be _I$}
Let $K\subset G$ be a nice (in the sense of \cite{BD}) compact open subgroup and
$$\mcH =\SS(G)^{K\times K},\ \  \mcH _L=\SS(L)^{ K\cap L\times K\cap L }$$
be the corresponding  Hecke algebras.
We use an arbitrarily chosen Haar measure on $G$ to identify $\mcH$ with
the space of compactly supported $K$ bi-invariant measures on $G$, thereby
endowing it with an algebra structure and the same applies to $\H_L$.

It is clear that
 $\SS (X_I)= i_I\circ r^- _I(\SS(G))$.
 We set  $M:=\SS (K\bc X_I/K)$.
We define  $T\in Aut(M)$ by $T=T^{(1,\dots,1)}$ (notations
of Definition \ref{T}).
Fix also a finite subset $\{ h_q\} ,q\in Q$ in $\mcH^{\ot 2}$ which generates $\mcH^{\otimes 2}$
 as a ring.

Let $C:=  \cup _{q\in Q}\supp(h_q) \subset G\times G$, 
 let $V\subset N_I$ be an open neighborhood of $S^0_I$ as in Proposition 
\ref{compatibility} and set
$$M^0:=\{ f\in M^K \ |\  \supp (f)\subset V\}.$$
Let $N=\mcH$, and $a:M^0\to N$ be the map coming from the admissible system $(\tau_I)$.

\begin{Prop}  The above data of $A=\mcH^{\otimes 2}$ modules $M$, $N$, the subspace
$M^0\subset M$ and automorphism $T$
 satisfy  conditions of Lemma \ref{ext}.
\end{Prop}

\proof  Noetherian property of $\mcH^{\ot 2}$ 
and finite generation of $M= \SS (K\bc X_I/K)$
follows from \cite{BD} (see in particular Remark 3.11 in {\em loc. cit.}).
The linear map $a:M^0\to N$ satisfies  compatibility condition 
of Definition \ref{compadef} when $h=h_q$
 by Proposition \ref{compatibility}. It follows that 
that it holds also when $h=h_{q_1}\cdots h_{q_n}$ for any $q_1,\dots, q_n\in Q$,
hence it holds for any $h\in \mcH^{\ot 2}$. 
 \medskip
 
Now Lemma \ref{ext} implies our main existence result:

\begin{Cor} \label{Bexists}
There exists  unique $\mcH^{\ot 2}$-covariant map $\Be ^K_I:M^K\to \mcH$ such that $\Be _I(T^{-n}f)=\tau_I
^* (T^{-n}f)$ for all
$f\in M^0$, $n\gg 0$.

\end{Cor}

\begin{Def}\label{B} We denote by $\Be _I:\SS(X_I)\to \SS(G)$ the linear operator whose restriction to $\SS (K\bc X_I/K)$
is equal to  $\Be ^K_I$ for all nice compact open subgroups $K\subset G$.\vspace{2mm}

\end{Def}

We will sometimes refer to $\Be _I$ as Bernstein's map.

\subsection{The induced map on $U_P$ coinvariants}
Fix a parabolic $P=P_I$. Consider the projection $X_I\to G/P_I\times G/P_I^-\overset{pr_2}{\To}
G/P_I^-$. Let $X_I^0$ be the preimage of the open $U_P$ orbit on $G/P_I^-$ under the composed map. Notice that the right action of $U_P$ on $X_I^0$ is free and $X_I^0/U_P\cong G/U_P$ canonically.
It follows that $\SS(X_I^0)_{U_P}=\SS(G/U_P)$.

Likewise, let  $^0X_I$ be the preimage of the open $U_P^-$
orbit under the composition $X_I\to G/P_I\times G/P_I^-\overset{pr_1}{\To}
G/P_I$. Then the $U_{P^-}$ action on $^0X_I$ is free and $^0X_I/U_P^-\cong G/U_{P^-}$ canonically,
hence $\SS(^0X_I)_{U_{P^-}}=\SS(U_{P^-}\bs G)$.

\begin{Prop}\label{composprop}
a)
The composed map
$$\SS(G/U_P)= \SS(X_I^0)_{U_P}\to \SS(X_I)_{U_P}\to \SS(G)_{U_P}=\SS(G/U_P),$$
where the first arrow is induced by the open embedding and the second one by the map $\Be _I$,
is equal to identity.

b)
The composed map
$$\SS(G/U_{P^-})= \SS(^0X_I)_{U_{P^-}}\to \SS(X_I)_{U_{P^-}}\to \SS(G)_{U_{P^-}}=
\SS( G/U_{P^-}),$$
where the first arrow is induced by the open embedding and the second one by the map $\Be _I$,
is equal to identity.
\end{Prop}

\proof Consider first the composed map in (a).
In view of Claim \ref{compfac} and Claim \ref{qUP} this map restricted to $K$-biinvariant
 functions on a neighborhood 
of $G/U_IZ_I\subset \overline{G/U_I}$ comes from an admissible map $\overline{G/U_I}\supset V\to
U\subset N_{\overline{G/U_I}}(G/U_IZ_I)$. However, it is easy to see that there exists
an algebraic isomorphism
$$N_{\overline{G/U_I}}(G/(U_IZ_I))\cong \overline{G/U_I}$$
such that the composed map
$$ X_I^0 \imbed N_{\Gb}(\Gb^0(I)) \to N_{\overline{G/U_I}}(G/(U_I Z_I))\cong \overline{G/U_I}$$
coincides with the projection $X_I^0\to G/U_I\subset \overline{G/U_I}$. 
(Here the second arrow in the displayed formula
 is the differential of the projection $\Gb^0(I)\to \overline{G/U_I}$).

  Thus identity map
$\SS(G/U)\to \SS(G/U)$ comes from an admissible map for
$\overline{G/U_I}$. It follows that the map in part (a) restricted to $K$-biinvariant functions on some neighborhood of the zero section $G/(U_IZ_I)$ equals identity. Since the map
is $G$ equivariant, it is equal to identity on all functions. This proves (a), part (b) is similar. \qed

\subsection{The map $\Be _I$ and $K_0$ cosets}
\label{53}

Let $K_0 \subset G$ be the standard open compact, i.e. $K_0$
is the group of $O$ points in the standard $O$ form of $G$
(recall that $G$ is assumed to be split).

Recall the canonical bijections $K_0\bs G/K_0=X_*(\mT)/W=X_*(\mT)^+, K_0\bs X_I/K_0=X_*(\mT)_{W_I}=X_*(\mT)^+_I$, where $W_I$ is the corresponding
parabolic Weyl group and $X_*(\mT)_I^+$ is the set of coweights positive on simple roots in the Levi. 
These are constructed by 
fixing a maximal torus  $T\subset G$ and a Borel subgroup $B$ containing $T$.
%
 Without loss of generality we can assume that the parabolic $P_I$, $I\subset \Sigma$ contains $B$.
Thus we get an embedding $T\to (G/U_I\times G/U_I^-)/L_I=X_I$. Let $\iota_I$ be
 the composition
of this embedding with the embedding $X_*(T)\to T$ sending  a cocharacter $\nu$ to $\nu (\pi)$,
where $\pi$ is a uniformizer.
For $\la\in X_*(\mT)^+$, $\nu\in X_*(\mT)_I^+$ the subsets $G_\la=K_0 \iota(\la) K_0$, 
$(X_I)_\nu=K_0 \iota_I(\nu) K_0$ do not depend on the auxiliary choices (here we abbreviated
$\iota=\iota_\emptyset:\Lambda \to G=X_\emptyset$).

Assume that $K$ is normalized by the maximal
compact subgroup $K_0$. Then $K_0\times K_0$ acts on $K\bs G_\la/K$.
It is easy to see that if $\la$ is such that  $(\la, \alpha)\gg 0$
for any root $\alpha$ in the radical of $P$, then 
the stabilizer of the point $K\iota(\la)K$ in $K_0^2$ 
 contains $K_0^+\times K_0^-$,
where $K_0^+=K_0\cap Rad(P_I)$, $K_0^-=K_0\cap Rad(P_I^-)$.
Since the group $K$  admits a triangular decomposition:
  $K=K_+ K_L K_-$ where $K_L=L\cap K$, $K_+=K\cap Rad(P_I)$,
 $K_-=K\cap Rad(P^-_I)$, it is easy to see that this
stabilizer coincides with the stabilizer of the corresponding point
in $K\bs (X_I)_\la/K$. Therefore
  there is a natural $K_0^2$-equivariant bijection
$\Psi_P:K\bs (X_I)_\la/K\isor K\bs G_\la /K$.

\begin{Lem} \label{cosets} Fix a congruence subgroup $K\subsetneq K_0$.
There exists $N>0$ such that for any $\lambda \in \Lambda^+$ satisfying 
$(\la,\alpha)>N$ for all roots $\alpha$ in the radical of $U_I\subset P_I$, the following holds:

a) The map $\Be _I$ sends
$\SS^{K\times K}((X_I)_\la)$ to $\SS^{K\times K}(G_\la)$ and induces an isomorphism
$\Be _I: \SS^{K\times K}((X_I)_\la) \isor \SS^{K\times K}(G_\la)$.

b) Moreover, the map $\Be _I|_{\SS(X_I)_\la}$ coincides with the map induced by 
bijection $\Psi_P$.

\end{Lem}

\proof For an integer  $N$ let $(X_I)_N$ be the union of $(X_I)_\la$
over $\la$ satisfying the condition of the Lemma.  These
 sets form a system of fundamental neighborhoods of the  closure $S_I$
of the corresponding stratum. For each stratum $S_J$ in the closure
 choose a neighborhood $V_J$ of $S_J$ in $G$. Without loss of generality we can assume that for each $J$ and an admissible system of maps
$\phi_J: V_J\to X_J$ the map $\Be _J|_{\SS(\phi_J(V_J))}$ coincides with $\phi_J^*$. We 
can also assume
that $V_J\subset V$ where $V$ is as in Proposition \ref{compatibility} for 
$C=K_0^2$.
Thus the maps $\phi_J^*$ on $K$-biinvariant functions are $K_0^2$ equivariant.

Fix $T\subset B$ as in section \ref{53}.
The choice of Borel subgroup $B$ defines an isomorphism
$T\isor \mT$. 
In view of Claim \ref{clXaff}a),
the resulting embedding $(\Gm)^\Sigma \to G$ extends to an embedding $
\Tb\to \Gb$.
Moreover, the intersection of its image with any stratum $G_I$ (notations
of section \ref{woncomp})  is 
 a single $T$-orbit.
  It is easy to see
that there is an admissible system for $\Tb$ 
where each map is a $T$-equivariant open embedding sending $\iota_I(\nu)$ to $\iota_J(\nu)$ for
$I\subset J$; here we use that the image of $\iota_I$ lies in 
$$N_{(\Aone)^\Sigma}((\Aone)^\Sigma
\cap \Gb_I)\subset N_{\Gb}(\Gb_I)\supset X_I.$$
 It also follows from the definitions that under the closed embedding this admissible system 
 is compatible with the one for $\Gb$. This yields
part (b) of the Lemma, which clearly implies part (a). \qed





\section{Second adjointness}\label{2ndadj}
\subsection{Basic notations}\label{61}
Notice that the spaces $X_I$, $G$ carry $G^2$ invariant measures (though $Y_I$ does not) 
and $G/U_I$ carries
a $G$-invariant measure. We fix such measures and identify the space $\SS$ of locally constant
compactly supported functions on $X_I$, $G$, $G/U_I$ with the space of locally constant
compactly supported measures.
 
We let $\MM(Y_I)$ be the space of locally constant compactly supported
 sections of the locally constant sheaf
$\mu_{pr_2}$ of fiberwise measures with respect to the second projection $pr_2:Y_I\to G/P_I$.
One can interpret elements
$\phi \in \MM(Y_I)$ as operators
$\hat \phi : \SS(G/U_I)\to \SS(G/U_I)$ where 

$$ \hat \phi (\psi )(\bar g):=
\int _ {\bar g'\in G/U_I}\phi (\bar g,\bar g') \psi (\bar g')dg',\psi \in \SS(G/U_I).$$

\begin{Def} We define the  {\it action map}
$$\A_I:\SS(G)\to \MM(Y_I),\ \ \  f\mapsto \pi_{2*} \pi_1^*(f),$$
where $\pi_1:(G\times G)/P_I\to G$, $\pi_1(g_1,g_2)=g_1g_2^{-1}$;
$\pi_2:(G\times G)/P_I\to Y_I$, $\pi_2((g_1,g_2)P_I)=(g_1,g_2)$
and $\pi_{2*}$ is the fiberwise integration map $\SS((G\times G)/P)\to \MM (Y_I)$;
here we use  the observation that the tensor product of the sheaf $\mu_{\pi_2}$ 
of fiberwise measures for the map $(G\times G)/P_I
\to Y_I$ and  $\pi_2^*(\mu_{pr_2})$ is canonically trivialized.
\end{Def}

\begin{Rem} 
 In terms of the interpretation of elements in $\MM(Y_I)$ as operators, 
 the map $\A$ is easily seen to correspond to the action of the Hecke algebra on the 
 universal principal series $\SS(G/U_I)$.
\end{Rem}
Let  $P=LU_P\subset G$ be a  parabolic  subgroup and $P^-=LU_P^-$ be the opposite parabolic  subgroup.
Let $i_I$, $r_I$ be, respectively, the normalized parabolic induction and the normalized Jacquet functors with
respect to $P_I$ and $i^-_I$, $r^-_I$ be those with respect to $P^-_I$.\vspace{2mm}

\subsection{The adjunction maps}
For a smooth $L_I$-module $N$ one defines canonical morphisms
$$can : r_Ii_I(N)\to  N,\, Can: N \to r^-_Ii_I(N).$$
To recall these observe that 
elements of the space of induced representations are $N$-valued functions on $G$, the map
$can$  comes from restriction of functions
to the closed subset
$P_I\subset G$, while $Can $ comes from push-forward of compactly supported
 functions from the open subset
$U_{P^-_I}P_I \subset G$.
 Using these canonical morphisms one can
define for any smooth $G$-module $M$ and a smooth $L_I$-module $N$ the following maps:
\begin{equation}\label{Phi}
Hom (M, i_I(N))\to Hom(r_I(M), N), \ \ \phi\mapsto can \circ r(\Phi)
\end{equation}
\begin{equation}\label{Psi}
 Hom (i_I(N), M )\to Hom(N, r_I(M^-), N),\ \
\psi \mapsto r^-_I(\psi)\circ Can.
\end{equation}

Frobenius  adjointness amounts to the  morphism \eqref{Phi}
 being an isomorphism.
This is a standard fact, one of possible proofs is as follows. It suffices to define
$\alpha: Id\to i_Ir_I$ so that the compositions $ i_I\to  i_I\circ r_I\circ i_I \to  i_I$, $r_I\to r_I\circ i_I\circ
r_I \to r_I$ are equal to identity.
Here the arrows $i_I\to  i_I\circ r_I\circ i_I $,
 $r_I\to r_I\circ i_I\circ r_I$
come from $\alpha$, and the arrows $ i_I\circ r_I\circ i_I \to  i_I$, $r_I\circ i_I\circ
r_I \to r_I$ come from \eqref{Phi}.

The morphism $\alpha$ comes from the map $\A_I:\SS(G)\to \MM(Y_I)$, since it is
 easy to see
 that $\MM(Y_I)\otimes _\H M=i_Ir_I(M)$, and it is obvious that $\SS(G)\otimes _\H M=M$.
 The compatibilites are easy to check.

 \subsection{Second adjointness}
We apply a similar strategy to show that
 \eqref{Psi} is also an isomorphism.

 We have  isomorphisms
$i_I\circ r^-_I(M)\cong \SS(X)\otimes_G M$ and $M\cong \SS(G)\otimes_G M$. Thus the map $\Be _I$ yields
for every $M\in Sm(G)$ a  map
$$\Be _I(M):i_I\circ r^-_I(M)\to M.$$
 In other words Bernstein's  morphism $\Be _I$ defines a
 morphism $\beta :i_I\circ r^-_I\to Id$.
Consider the compositions:

$$\nu _I:i_I\overset{Id\otimes Can}{\To} i_I\circ r^-_I\circ i_I\overset{\beta \otimes Id}{\To}  i_I$$
and

$$\tau _I: r^-_I\overset{Can \otimes Id }{\To} r^-_I\circ i_I\circ {r^-}_I\overset{Id \otimes \beta}{\To}
{r^-}_I$$

\begin{Thm}\label {adj}
a) The  morphisms $\nu _I$ and $\tau _I$ are the identity morphisms.

b) The map \eqref{Psi} is an isomorphism. In particular, $i_I$ is the left adjoint of $r_I^-$.

\end{Thm}
\proof (b) follows from (a) by a standard argument. To check (a) observe that the maps
of functors $\nu_I$, $\tau_I$ come from the maps of bimodules considered, respectively,
in parts (a) and (b) of Proposition \ref{composprop}. Thus part (a) of the Theorem follows from that Proposition. \qed

\section{Bernstein's map and intertwining operators}\label{intop}

In this section we formulate some properties of the maps introduced above and state a result which says that the maps $\A$ and $\Be $ are related by the intertwining operator. The proof occupies the next
 two sections.
 
 To simplify notation we only treat the case of the spaces $X,Y$ attached to the Borel subgroup,
 we expect that similar results can be proved in a similar way for spaces $X_I$, $Y_I$. 

Let $\SS'(G)$, $\SS'(X)$ denote the space of functions which are invariant under $K\times K$
for some open compact subgroup $K$ (but are not necessarily compactly supported).
Our choice of Haar measure defines a pairing between $\SS(X)$ (respectively,
$\SS(G)$) and $\SS'(X)$ (respectively, $\SS'(G)$), identifying $\SS'(X)$, $\SS'(G)$ with smooth dual of $\SS(X)$, $\SS(G)$ respectively. We let $\Be ^\star:\SS'(G)
\to \SS'(X)$ denote the operator dual to $\Be $.

\subsection{Bounded supports}



 Recall that $\Xb_{\aff}$ denotes (the space of $F$-points of) 
 the affine closure of the quasi-affine algebraic variety $X$.
Let us  say that a subset $C\subset X$   
 is bounded if
the closure of $C$ in $\Xb_{\aff}$ is compact.

Let $\SS_b(X)$
denote the space of locally constant functions on $X$ 
 which are supported on a bounded subset in $X$.

\begin{Prop}\label{Bstim}
For $f\in \H$ the support of $\Be 
^\star(f)$ is bounded, i.e.

 $\Be ^\star: \SS(G)\to \SS_b(X)$.
\end{Prop}

In view of potential applications we mention a more precise version of this result
in a special case; it will not be used in this paper.

\begin{Prop}\label{BstimDr}
The support of $\Be 
^\star(\delta_{K_0})$ is contained in $\Xb_{\aff}(O)$.
\end{Prop}

The proof of Proposition \ref{BstimDr} appears in section \ref{Plf}.

\subsection{Intertwining operators}\label{secin}
Fix $w\in W$.
Recall the Radon correspondence $\fC_w$ appearing in Lemma \ref{fCwdef}.
We define the intertwining operator (or Radon transform) $I_w:\SS(X)\to \MM'(Y)$, $I_w:f\mapsto pr_{2*}^w (pr_1^w)^*(f)$ where $pr_1^w:\fC_w\to X$, $pr_2^w:\fC_w\to Y$ are the projections
and $\MM'(Y)$ is the space of locally constant sections of the sheaf $\mu_{pr_2}$ (see
\S \ref{61}). 

\begin{Prop}\label{Idefined}
Let $S\subset X$ be a closed bounded subset. Then the map $pr_2^w:(pr_1^w)^{-1}(S)\to Y$ is proper.
\end{Prop}

\begin{Cor}\label{IwCb}
The intertwining operator $I_w$ naturally extends to a map $I_w:\SS_b(X)\to \MM'(Y)$.

\end{Cor}

Let us say that a subset in $Y$ is $w$-bounded if its closure in $\Yb_w$ is compact.

Let $\MM_b^w(Y)$ denote the space of locally constant sections of  the sheaf $\mu_{pr_2}$
 whose support is $w$-bounded.

\begin{Prop}\label{Iisom}
a) The map $I_w$ sends $\SS_b(X)$ to $\MM_b^w(Y)$.

b) The map $I_w:\SS_b(X)\to \MM_b^w(Y)$ is an isomorphism.
\end{Prop}

\subsection{Main result}
\begin{Thm}\label{AIBst}
For any $w\in W$ we have $\A=I_w \Be ^\st$.
\end{Thm}

\begin{Rem} The latter equality (in the special case $w=1$)
 resembles the result of \cite[Corollary 6.2]{BFO}.
More precisely, in  {\em loc. cit.} one finds an isomorphism of two functors between the categories of $D$-modules. Following a standard analogy between maps of function spaces and functors
on the categories of $D$-modules given by "the same" correspondence one gets that
one of the two functors considered in {\em loc. cit.} is analogous to the map $I_w^{-1}\A$. 
The other functor in \cite{BFO} is that of {\em nearby cycles,} or specialization.
The characterization of $\Be $ via co-specialization of functions on an $F$-manifold to normal cone
 makes it natural to consider $\Be ^\st$ as an analogue of specialization functor between the categories
 of $D$-modules. 
 It would be interesting to find a precise mathematical statement underlying these heuristic 
 considerations.
\end{Rem}

The proof of the Theorem is given in  section \ref{prThm}.

\begin{Cor}\label{Bform}
For any $w\in W$ we have $\Be =\A^\st (I_w^\st)^{-1} $.
\end{Cor}
\section{Some properties of Radon correspondence}\label{proRad}



\subsection{Proof of Propositions \ref{Idefined}, \ref{Iisom}}\label{proofIisomb}

 Proposition  \ref{Idefined} follows from Lemma \ref{imageof}.

Proposition \ref{Iisom}(a) follows from Proposition \ref{prprop}.

In order to prove 
Proposition \ref{Iisom}(b) we show the following more precise result.

Along with $I_w$ we will also consider the adjoint operator $I_w'=pr^w_{1*}(pr^w_2)^*:\MM_b(Y)\to \SS'(X)$.

\begin{Prop}\label{sigmasi}
There exists an element $\sigma$
in the group algebra of the torus $T$ with the following properties.

i) The element $\sigma$ 
is a  product of elements of the form
$[\alpha_i]-c_i$ where $c_i$ is a constant and $[\alpha_i]\in T$ is a representative
of the coset of $T^0$ corresponding to a positive coroot.

ii) For $f\in \SS(X)$ the element $I_w(\sigma
(f))$ has compact support.

iii) For $f\in \SS(X)$ the element $I_w'I_w(\sigma
(f))$ (which is well defined by ii)
has compact support.

iv) For a character $\chi$ of $T^0$ let $\SS_\chi(X)\subset \SS(X)$
be the subspace of $\chi$ semi-invariants with respect to the action of $T^0$.
 For every $\chi$ the restriction of the map $f\mapsto I_w'I_w(\sigma(f))$ to $\SS_\chi$
 equals the action of 
an element $\sigma_\chi$ in the group algebra of $T$, where $\sigma_\chi$
is a product of nonzero elements of the form $b_{\chi,i} [\alpha_i]-c_{\chi,i}$ where
$\alpha_i$ is as in (i). 
\end{Prop}

\proof Fix a minimal decomposition for elements $w$, $w^{-1}w_0$. This defines a
presentation of $I_w$ as a composition of $\ell(w)$ simpler correspondences.
Thus it suffices to prove a similar statement for each of these simpler correspondences.

This reduces to the following well known properties
 of Radon transform. 
  
 \begin{Claim}\cite[\S II.2.5, II.2.6]{GGPS} \label{Radpl}
 Let $R$ denote Radon transform for functions on the plane,
$R(f)(y)=\int f(x_0+ty)dt$ where $x_0$ is such that $\langle x,y\rangle =1$
for a fixed skew-symmetric bilinear pairing $\langle\ ,\ \rangle$.
 
 a) If $f\in \SS(F^2\setminus \{0\})$ is such that the integral of $f$ over any line 
 passing through zero vanishes, then $R(f)$ has compact support.
 
 b) If $\pi\in F$ is the uniformizer and $q$ is the cardinality of the residue field
 of $F$, then $\pi(f)-q^{-1}f$ satisfies the assumption of (a) for any $f\in \SS(F^2\setminus \{0\})$, where for $t\in F^\times$ we write $t(f)(x):=f(t^{-1}x)$.
 
 c) Let $\Phi$ denote Fourier transform for functions on the plane.
 Then there exists  a rational scalar valued function $\gamma$ on the set of multiplicative
 characters, such that for any $f\in  \SS(F^2\setminus \{0\})$ satisfying the assumption of (a) and 
 any multiplicative character $\chi$ in the domain of definition of $\gamma$
 we have $$\gamma(\chi) \cdot  \overline{R(f)}_\chi = \overline{\Phi(f)}_\chi.$$ Here
 for $\phi\in   \SS(F^2\setminus \{0\})$ we let $\overline{\phi}_\chi$ denote 
 the image
 of $\phi$ in the space $ \SS(F^2\setminus \{0\})_\chi$ of $\chi$-coinvariants with respect to the dilation action of $F^\times$.
 
 We also have
  $\Phi^2=Id$.
 \end{Claim}
 
\begin{Rem}
The above proof does not apply in the case when Borel subgroup is replaced by a parabolic one.
It is possible to give an alternative proof admitting such a generalization.
\end{Rem}

\subsubsection{Proof of Proposition \ref{Iisom}(b)}
In view of Proposition \ref{sigmasi}  $I_w$ induces
an isomorphism between the spaces $\SS(X)$ and $\MM(Y)$ tensored with
 localization of the group algebra of the torus $T$ by some elements
of the form $[\alpha_i]-c_i$ where $\alpha_i$, $c_i$ are as in Proposition
\ref{sigmasi}. The subset $\{\alpha_i^n\ |\ n\geq 0\}\subset T$ is 
bounded, which shows that for $f\in \SS_b(X)$ the infinite sum
$\sum\limits_{n=0}^\infty c_i^{-n}\alpha_i^n(f)$ is a well defined
element of $\SS_b(X)$. Thus $[\alpha_i]-c_i$ acts by an invertible 
operator on $\SS_b(X)$, a similar argument shows that it also acts
by an invertible operator on $\MM_b^w(Y)$. Proposition \ref{Iisom}(b) follows.

\medskip

The rest of the section is devoted to the proof of Proposition \ref{Bstim}.

\subsection{The Bernstein center and supports}
\subsubsection{The Bernstein center}
Let $\cZ$ denote the Bernstein center of $G$. 
 Recall \cite{BD} that $Spec(\cZ)$ is the union of connected components,
and each component has the form $Spec(\cZ_k)=T_k/{W_k}$  where the torus $T_{L_k}$ is  dual to $L_k/L_k'$ for   a Levi subgroup $L_k$  in $G$, 
while $W_k$ is a finite group acting on $T_k$ (in fact, $W_k$ is a subgroup in $T_{L_k}\ltimes W_k$ where $W_k$ is the Weyl group; here we used the standard notation $L'=[L,L]$
for the commutant). Thus for every $k$, the summand $\cZ_k$
is a subalgebra in $\Ce[X_*(L_k/L_k')]$.
We have  embeddings
$X_*(\A_k)\to X_*(L_k/L_k')$, $   X_*(\A_k)\to X_*(\mT)$
 where $\A_k$ is the center of $L_k$; the first embedding has a finite index.
 Thus $X_*(L_k/L_k')\subset X_*(\mT)_\Qu:=X_*(\mT)\otimes \Qu$ canonically.

We will need the following property of
elements in the Bernstein center $\cZ$ which follows directly from the 
description of $\cZ$ in \cite{BD}.

\begin{Claim}\label{centcharn_a}
 For any element $h\in \cZ$ there exists an element
$h_L\in \cZ(L)$ such that the left action of $h$ on $\SS(G/U_P)$ coincides with the action of $h_L$
coming from the right action of $L$. 
\end{Claim}

\subsubsection{Filtrations by support}
We define an increasing filtration on
$\cZ$ as follows.
 For $\mu\in X_*(\mT)^+$ set \vspace{2mm}

$X_*(\mT)_\Qu ^{\leq \mu}=\{\la\in X_*(\mT)_\Qu \ |\ (w(\la), \alpha)\leq (\mu,\alpha) $ for all simple coroots $\alpha,\ w\in W\}$,\vspace{2mm}

 and define

 $$(\cZ_k)_{\leq \mu}=\cZ_k \cap \Ce[X_*(\mT)_\Qu ^{\leq \mu}],\ \ \
 \cZ_{\leq \mu} =\prod (\cZ_k)_{\leq \mu}.$$
 This equips $\cZ_k$, $\cZ$ with filtrations indexed by the partially ordered semi-group of dominant weights.
It is easy to see that the Rees ring $\oplusl_{\mu}  (\cZ_k)_{\leq \mu}$  is  finitely generated for each $k$.

\medskip

We will need an auxiliary result relating this filtration to a filtration on the Hecke algebra.

Recall notations $K_0$, $G_\la$ etc. from section \ref{53}.

We let $G_{\leq \la}=\bigcup \limits _{\mu\leq \la} G_\mu$ and we let
$\H_{\leq \la}=\{h\in \H \ |\ \supp(h)\subset G_{\leq \la}\}$; here $\leq$ denotes
the standard order on dominant coweights.

Fix a congruence subgroup $K\subsetneq K_0$. Thus $K$ is a normal open subgroup in 
$K_0$ and $K$ is nice in the sense of \cite{BD}.
We change notation and write 
$\H$ instead of $\H_K$ and $\H_{\leq \la}$ instead of $\H_{\leq \la}\cap \H_K$.

\begin{Prop}\label{comp} 
 There exists $\la_0$ such that for $\la\geq \la_0$ and any $\mu$ we have

$$\cZ_{\leq \mu}\cdot \H_{\leq \la}\subset \H_{\leq \la+\mu}.$$
\end{Prop}

The proof of the Proposition is preceded by two auxiliary statements.

The first one  is a property of the element
$h_L$ introduced in Claim \ref{centcharn_a}, it is an immediate consequence of characterization
of $h_L$ formulated in that Claim.

\begin{Claim}\label{centcharn_b}
For $h\in \cZ_{\leq \mu}$ let $h_L$ be as in Claim \ref{centcharn_a}.
Then the support $\supp(h_L)$
satisfies the following condition. Let $L_c\subset L$ be the subgroup generated by compact subgroups, so that
$L/L_c\cong X_*(L/L')\subset X_*(\mT)_\Qu$. The image of $\supp (h_L)$ in $L/L_c$ is contained in
$ X_*(L/L')_{\leq \mu}$. 
\end{Claim} 

To state another Lemma, recall that
we can identify the double quotient $K_0\bs X_I/K_0$ with
$X_*(\mT) /W_I=X_*(\mT)_I^+$, where $W_I$ is the corresponding
parabolic Weyl group and $X_*(\mT)_I^+$ is the set of coweights positive on simple roots in the Levi subgroup $L_I$.
For  $\nu\in X_*(\mT)_I^+$ we denote by  $(X_I)_\nu$  the corresponding $K_0^2$ orbit.\vspace{2mm}


\begin{Lem} \label{ZXP}
Fix $I\subset \Sigma$.
For any $\mu \in X_*(\mT)^+$ there exists a  finite set $S_\mu$
of linear combinations of the form $\sum r_i\alpha ^\vee _i,\, i\in I$ such that
$\cZ_{\leq \mu}\SS(X_I)_\nu^{K\times K} \subset \sum \SS(X_I )
_\eta, \ \ \eta\in \nu + {\mathrm{conv}}
(W(\mu)) + S_\mu$ where ${\mathrm{conv}}$ is the convex hull.
\end{Lem}

\proof We have a fibration $X_I 
\to G/P_I^-$ with fibers $G/U_I$; the set $K_0\bs G/(U_I\cdot K_0^{L_I})$ maps isomorphically to $K_0\bs X_I/K_0$. The left action of $\H$ we consider here is fiberwise, thus it suffices to prove a similar statement about the action on $\SS(G/U_I)$.

An element $h\in \H(L)$ supported on the preimage of a given $\zeta\in
 X_*(L/L')$ sends $\SS(G/U_I)_\nu$ to the sum of $\SS(G/U_I)_{\eta}$ where $\eta\in \nu+\zeta+S$
 for some finite subset $S=S(h)$ of linear combinations of $\alpha ^\vee _i$, $i\in I$. Thus
  Lemma follows from Claim \ref{centcharn_b}. \qed

\subsubsection{Proof of Proposition \ref{comp}}
Since the  Rees ring $\oplusl_{\mu} ( \cZ_k)
_{\leq \mu}$ is  finitely generated for each $k$ and ${\cZ_k} \H =0$ for almost all $k$,
it suffices to show that for fixed $k$, $\mu$ and sufficiently large $\lambda$ we have
\begin{equation}\label{cZk}
(\cZ_k)_{\leq \mu}\cdot \H_{\leq \la}\subset \H_{\leq \la+\mu}.
\end{equation}

Fix $N\in \Zet_{>0}$ as in Lemma \ref{cosets}.
We can and will assume that $N>2r_i$ for $r_i$ as in Lemma \ref{ZXP}, $\mu\in M$.
We can and will assume also that $\< \lambda, \alpha\> >2N$ for any simple root $\alpha$.\vspace{2mm}


Fix $\nu \in X_*(\mT)$, $\nu\leq \la$ and let $f$ be a $K$-biinvariant function 
supported on $G_\nu$. Let $P=P_I$ be the parabolic subgroup such that the simple roots in its Levi are exactly those simple roots $\alpha$ for which $(\alpha, \nu)\leq N$.
By Lemma \ref{cosets}(a) we have $f=\Be _I(f_P)$ for some $f_P\in
\SS(X_I)_\nu$.
Thus for $z\in Z_{\leq \mu}$ we have:
$z(f)=\Be _I(zf_P)$. In view of  Lemma \ref{ZXP}, $z(f_P)\in \sum \SS(X_I)_\eta$ with
$\eta\in \nu + {\mathrm{conv}}(W(\mu))+ S_\mu$ for   a fixed finite set $S_\mu$
of linear combination $\sum r_i \alpha_i$ where $\alpha_i$ are coroots of the Levi.

To finish the proof of Proposition \ref{comp} it suffices to check that 
$\nu +\sum r_i\alpha_i\leq \lambda$.
To do this, notice that  for every simple coroot $\alpha_i$ of
the Levi and the corresponding fundamental weights $\omega_i$ we have:

$$(\lambda - \nu, \omega_i)>\frac{1}{2}(\la-\nu, \alpha_i)> \frac{1}{2}(2N-N) =\frac{1}{2}N>r_i, $$ where the first inequality follows from the fact that $\la-\nu$ is a sum of positive coroots, while the other ones follow from the assumptions on $N$. 
 \qed\vspace{2mm}

In the next subsection we will need the following consequence of 
Proposition \ref{comp}.\vspace{2mm}

\begin{Lem}\label{pronos}
Fix an open compact subgroup $K\subset K_0$. There exists
a finite subset $S\subset X_*(\mT)$  such that for any $\mu\in X_*(\mT)_I^+$, $\phi\in \SS(X_I)_ \mu$ and
$\lambda \in X_*(\mT)^+$ and any $w\in W$
there exist $\psi_i\in \SS(X_I)$, $h_i\in \H$, such that

$\phi=\sum h_i(\psi_i)$, $supp(h_i)\subset G_{\leq \la}$, $supp(\psi_i)\subset (X_I)_{\mu-w(\la)+S}$.
\end{Lem}

\proof 
The space $\SS(X)$ is acted upon by the maximal torus $\mT$; for our current purposes it suffices to consider the action of the subgroup $X_*(\mT)\subset \mT$ (where the
embedding depends on the choice of a uniformizer $\pi\in F$).\vspace{2mm}

It follows from the description of $\cZ$ in \cite{BD}  that the action of an element
 $s\in \Ce[X_*(\mT)]^W$ on $\SS(X)^{K\times K}$ coincides with an action of some element in $s'\in \cZ$. Moreover, if
$s\in \Ce[X_*(\mT)]_{\leq \mu}$ for some $\mu \in X_*(\mT)^+$, then $s'$ can be chosen in $\cZ_{\leq \la}$.\vspace{2mm}

It is clear that $\SS(X)_{0}$ generates $\SS(X)$ as a module over $X_*(\mT)$. Choose a finite dimensional space of generators for $\Ce[X_*(\mT)]$ over $\Ce[X_*(\mT)]^W$, then applying it to $\SS(X)_{0}^{K\times K}$ we get a space of generators for $\SS(X)^{K\times K}$ over $\Ce[X_*(\mT)]^W$, let us denote it by $V$. We can also assume without loss of generality that $\Ce[X_*(\mT)]_{\leq \la} \SS(X)_0\subset \Ce[X_*(\mT)]^W_{\leq \la} V$ for all $\la\in X_*(\mT)^+$.\vspace{2mm}

Let us now choose $\la_0$ as in \ref{comp} and  choose $S$  such that
$\SS(X)_S \supset \Ce[X_*(\mT)]^W_{\leq \la_0}\cdot V$. We claim that this $S$ satisfies the condition of Lemma \ref{pronos}.
 Indeed, for $\la \geq \la_0$ we have

$$\SS(X)_{\leq \la} \subset \Ce[X_*(\mT)]_{\leq \la} V
=\Ce[X_*(\mT)]^W_{\leq \la-\la_0} \Ce[X_*(\mT)]^W_{\la_0} V\subset \H_{\leq \la}\SS(X)_S.$$

Here in the last inclusion we used the inclusion  $\Ce[X_*(\mT)]^W_{\leq \la-\la_0}\delta_K
\subset \cZ_{\leq \la-\la_0}\delta_K\subset \H_{\leq \la}$ which follows from \ref{comp}. \qed\vspace{2mm}

The last ingredient needed in the proof  of \ref{Bstim} is the following easy statement.

\begin{Claim}\label{compcl}
For any $\la_0\in X_*(\mT)$ 
the subset $X_{\leq \la_0}:=\cup_{\la\leq \la_0} X_{I,\la}$
has compact closure in $\Xb_{\aff}$.
\end{Claim}

\subsection{Completion of the proof of Proposition \ref{Bstim}}
Since the space of bounded distributions is clearly invariant under the $G\times G$ action, it suffices to prove the statement for $f=\delta_K$.


In view of Claim \ref{compcl} it is enough to show existence of $\la_0 \in X_*(\mT)_I^+$ such that for $\phi\in \SS(X)^{K\times K}$ the condition
$supp(\phi)\subset X_{I,\mu}$, $\mu \not \leq \la_0$ implies that $\Be (\phi)|_1=0$.

Assume $\phi\in \SS(X)_\mu$ and write
$\mu=\la_1-\la_2$, where $\la_1,\la_2\in \La^+$. Apply Lemma \ref{pronos} with $w=w_0$. We get that  $\phi=h(\psi)$, for some $h\in \H_{\leq \la_2'}$ and $supp(\psi)\subset X_{\la_1+S}$; here $\la_2'=-w_0(\la_2)$ is the dual weight.
We can and will assume without loss of generality that
the pairing of each coweight in $\la_1+S$ with any fundamental weight is larger than $N$ where
$N$ is as in Lemma \ref{cosets}.
This Lemma implies then that  $\supp(\Be (\psi))\subset G_{\la_1+S}$.
Thus $f*\Be (\psi)|_1=0$ for $f\in \H_{\leq \nu}$ unless $\la_1+\sigma \leq \nu'$ for some
$\sigma \in S$, $\nu'=-w_0(\nu)$. If we choose $\la_0$ so that $\la_0\not \leq -\sigma$
for $\sigma \in S$, we get that $\la_2=\la_1-\mu \not \geq \la_1 +\sigma$, so $h*\Be (\psi)|_1=0$ since $h\in \H_{\leq \la_2'}$. \qed

\medskip

\section{Proof of 
Theorem \ref{AIBst}}\label{prThm}

Recall that $\SS_b(X)$, $\SS_b(T)$
denotes the space of locally constant functions on $X$ (respectively, $T$) with bounded support;
$X_{w_1,w_2}$ denotes the corresponding $B\times B$ orbit on $X$.

Recall the 
 Jacquet and induction functors $r_I=r_{P_I}$, $i_I=i_{P_I}$ introduced
in section \ref{61}; when $P_I=B$ is a Borel we abbreviate this to $r$, $i$.
We will also abbreviate $r(\SS(G))$ to $r(G)$ etc.
 


\begin{Prop} The stratification of $X$ by $B\times B$ orbits $X_{w_1,w_2}$ induces
a filtration on $r(\SS_b(X))$ which splits canonically, yielding a canonical isomorphism:

 \begin{equation}\label{WW}
 r(\SS_b(X))
 \cong \oplusl_{W
 \times W
 } \SS_b(T).
\end{equation}

\end{Prop}

We denote the summand corresponding to $(w_1,w_2)\in W\times W$ by 
$r(\SS_b(X))_{w_1,w_2}$.

The Proposition  follows from the next Lemma which shows that for every 
$B\times B$ orbit $X_{w_1,w_2}\subset X$ the map
$r(\SS_b(X))\to \SS_b(T)$, $f\mapsto (pr_{w_1,w_2})_*(f|_{X_{w_1,w_2}})$ is well
defined; here $pr_{w_1,w_2}:X_{w_1,w_2}\to T\cong X_{w_1,w_2}/(U_B\times U_B)$
is the projection. It is clear that the direct sum of these maps provides
an isomorphism \eqref{WW} which splits the filtration. \qed

\begin{Lem}\label{bet}
 Let $\beth\subset X$ be a bounded subset. Fix $w_1,w_2\in W$ 
 and let
$X_{w_1,w_2}$ be the corresponding $B\times B$-orbit. Then the map
$pr_{w_1,w_2}:\beth\cap X_{w_1,w_2}\to T=X_{w_1,w_2}/U_B^2$ is proper and its image is bounded.

For every $w_1$, $w_2\in W$, every bounded subset in $T$ is the image of $\beth\cap X_{w_1,w_2}$ under $pr_{w_1,w_2}$
for some bounded subset $\beth\subset X$.



\end{Lem}

\proof Without loss of generality we can assume that $\beth$ is invariant under an open compact subgroup
in $G\times G$. Then properness  amounts to compactness of the intersection of $\beth$ with every
fiber of the projection $X_{w_1,w_2}\cap \beth\to T$, which  follows from the fact that the fibers of projection
$X_{w_1,w_2}\to X_{w_1,w_2}/U^2$ are closed in $\Xb_{\aff}$. The latter is a consequence
of Kostant--Rosenlicht Theorem 
 \cite[Proposition 2.4.14]{Sprb}, since each fiber is an orbit of an affine algebraic group
 acting on an affine algebraic variety.

The fact that the image of the projection is bounded  and the last statement follow respectively from parts b) and a) of Claim \ref{clXaff}. \qed

\subsection{The subquotient maps for $\Be ^\st$}

In the next statement we use identifications of $T$-torsors $U_B \bs G_w /U_B$ and $U_B\bs X_{w_1,w_2}/U_B$
where $w_1w_2w_0=w\in W$. Such an identification follows from Corollary
\ref{BrclCor}.

\begin{Lem}\label{compBst}
Suppose that $w=w_1w_2w_0$ and $\ell (w)+\ell(w_0)=\ell(w_1)+\ell(w_2)$.
Then the composition  $r(G)_{\geq w}\overset{r(\Be ^\st)}{\To}  r(\SS_b(X))\to
r(\SS_b(X))_{w_1,w_2}$ factors through a map
$r(G)_w\to r(\SS_b(X))_{w_1,w_2}$. This map equals the canonical embedding $\SS(T)\to \SS_b(T)$.
\end{Lem}

\proof The composed map in question
factors through a map
$r(G)_w\to r(\SS_b(X))_{w_1,w_2}$ because there are no nonzero $T^2$ equivariant
maps $(r(G)_{<w})\to r(\SS_b(X))_{w_1,w_2}$.


The statement  readily follows from the following formula:
\begin{equation}\label{prw12}
pr_{w_1,w_2*}(\Be ^\st(f)|_{X_{w_1,w_2}})=\overline{f};
\end{equation}
here $f\in \SS(G)$ is such that $f|_{G_{<w}}=0$ and $\overline{f}$ is the image of
$f$ in $r(\SS(G))_w=\SS(T)$. Notice that the direct image is well defined since
the map
$$\supp (\Be ^\st(f))\cap X_{w_1,w_2}\to T=X_{w_1,w_2}/U_B^2$$
 is proper in view
of Lemma \ref{bet} and Proposition \ref{Bstim}.

Since the map $pr_{w_1,w_2*}$ factors through coinvariants with respect to the action of
$U_B\times U_B$, we have $pr_{w_1,w_2*}(\Be ^\st(f)|_{X_{w_1,w_2}} )= F(\overline{f})$ for some
$T\times T$ equivariant map $F:\SS(T)\to \SS_b(T)$. Thus \eqref{prw12} would
follow if we show that restrictions of both sides to a non-empty open subset $C\subset T$ (which may depend on $K$ but not on $f$) coincide. The construction of the map $\Be $ in section \ref{consB}
(see 
Corollary \ref{Bexists}) makes it clear that for some neighborhood
$V_G$ of $Z=\BB\times \BB$ in $\Gb$, a neighborhood $V_X$ of
the zero section in the normal bundle $N_{\Gb}(Z)$
and an admissible bijection $\tau:V_X\to V_G$ we have $\Be (f)=\tau_*(f)$ for any $f\in \SS(X)^{K\times K}$, $\supp(f)\subset V_X$. In view of Claim \ref{compfac}, we will be done if we show that for some
$U_B\times U_B$ invariant subset $V\subset X_{w_1,w_2}$ we have $V\subset V_X$.
This follows from the next geometric Lemma.
\qed

\begin{Lem}
For a representation $M$ of $G$ let $\rho_M:  
X \to End(M)$ be the canonical map
 as in Claim \ref{clinEnd}.

a) Suppose that a subset $\fZ\subset X$ is such that for any $M$ the closure of
$\rho_M(\fZ)$ does not contain zero. Then for any
 neighborhood $V$ of the zero section  in $N_{\Gb}(Z)\supset X$ we have
$z_I^N(\fZ)\subset V$ for $N\gg 0$ (where we use notation of Definition \ref{T}). 

b) Every orbit of the group $U_B \times U_B\times T_0$, where $T_0$ is a compact subgroup
 in the abstract Cartan, satisfies the conditions of part (a).
\end{Lem}

\proof  (a) follows from Claim \ref{clinEnd}, in view of the following easy observation:
 given finite dimensional vector spaces $\VV_i$ over $F$, and a subset $\fZ$ in $\prod \VV_i^0$
such that the image of $\fZ$ under the $i$-th projection does not contain zero in its closure,
the closure of $\fZ$ in $\prod \bbPt(\VV_i)^0$ is compact
(where we used notations of \ref{embinos}); hence for any neighborhood $V$ of
$\prod \bbP(\VV_i)$, multiplication by $\pi^{-N}$ sends $\fZ$ to $V$ for some $N$.

To prove (b)
 it is enough to show the same statement for a $U_B\times U_B$ orbit. The image
of such an orbit in $End(M)_i$ is also a $U_B\times U_B$ orbit different from $\{0\}$. Since an orbit
of a unipotent group on an affine algebraic  variety is Zariski closed by Kostant-Rosenlicht Theorem, we get the statement. \qed

\begin{Rem}
The implication in part (a) of the Lemma  is actually an "if and only if" statement, we only proved
the direction we need to save space.
\end{Rem}

The last auxiliary fact needed in the proof of  Theorem  \ref{AIBst} is the following  property of intertwining operators which is immediate from its definition.

For $w_1,w_2\in W$ we let  $r(\MM_b^v(Y))_{w_1,w_2}$ denote the corresponding subquotient
of the filtration on $r(\MM_b^v(Y))$ induced by the stratification of $Y$ by $B\times B$ orbits.

\begin{Claim}\label{compI}   The component  
$r(\SS_b(X))_{w_1,w_2}\to r(\MM_b^v(Y))_{w_1',w_2'}$ of $I_v$ equals identity if $w_1'=vw_1$, $w_2'=w_2 w_0 v^{-1}$ and $\ell(w_1)=\ell(v) +\ell(w_1')$,
  $\ell(w_2)=\ell(w_2')+\ell( w_0 v^{-1})$.


\end{Claim}


\subsection{Proof of Theorem \ref{AIBst}}

Comparing Proposition \ref{Bstim} with Corollary \ref{IwCb} we see that the composition
$I_w \Be ^\st$ is well defined.

Then Lemma \ref{compBst} together with Claim \ref{compI} show that $r(I_wB^\st)$ induces a quotient map
$r(\SS(G))_1 \to r(\MM_b^w(Y))_{1,1}$ which is equal to the canonical embedding
$\SS(T)\to \SS_b^w(T)$. By Frobenius adjointness  this implies that $I_wB^\st$ is the composition of $\A$ with the embedding $\MM(Y)\to \MM_b^w(Y)$. \qed

\section{Plancherel type formulas}\label{Plf}
In this last section we present an application of 
Corollary \ref{Bform}.

\subsection{The map $\Be $ and the Plancherel functional}\label{101}
The equality $\Be = \A^\st (I_w^\st)^{-1}$ allows one to write down a spectral expression 
for the value $f(1)$ for a function $f\in \SS(G)$ presented in the form $f=\Be (\phi)$, $\phi\in \SS(X)$.
To spell this out we need some preliminaries. 

Let $Char$ denote the space of characters of the torus $\mT$ viewed as an affine
ind-algebraic variety
over $\Ce$. We have $Char\cong \bigsqcup\limits_{\chi \in \mT(O)^\vee} Char_\chi$, where
$\mT(O)^\vee$
is the discrete set of characters of $\mT(O)$
and $Char_\chi$ is the space of characters whose restriction to $T(O)$ is equal to $\chi$.
Then $Char_\chi$ is a principal homogeneous space for  the dual torus
 $\mT^\vee$
over $\Ce$. Let $\OO^f(Char)=\bigoplus\limits_\chi \OO(Char_\chi)$ denote the (non-unital) ring of regular functions on $Char$ which vanish
on all but a finite number of components. We have a canonical isomorphism $\OO^f(Char)\cong \SS(\mT)$,
and  a noncanonical isomorphism $\OO^f(Char)\cong \bigoplus\limits_{\chi\in \mT(O)^\vee} \OO( \mT^\vee)$.

Recall the partial compactification $\mTb_{\aff}$ of the torus $\mT$ given by $\mTb=Spec(F[X^*(\mT)^+))$
where $X^*(\mT)^+\subset X^*(\mT)$ is the semigroup consisting of positive rational 
combinations of simple roots in $X^*(\mT)$. Let $\SS_b(\mT)$ be the space of distributions 
on $\mTb_{\aff}$ invariant with respect to an open subgroup in $\mT$  whose support is contained in a compact subset of $\mTb_{\aff}$. 
To give a spectral description of this ring introduce the affine
toric variety $\mTb^\vee=Spec(\Ce[X_*(\mT)^+])$, where 
$X_*(\mT)^+$ 
 is the dual cone to $X^*(\mT)^+$, i.e. 
$X^*(\mT)^+\subset X_*(\mT)=X^*(\mT^\vee)$ is the set of weights which are positive rational
linear combinations 
 of dominant weights of $\mT^\vee$. This defines a partial compactification
$\overline{Char}_\chi=(\mTb\times Char_\chi)/\mT$ of the principal homogeneous space $Char_\chi$. Let $\Ohat(Char_\chi)$ be the ring of functions on the punctured formal neighborhood of the divisor
$\partial Char_\chi=\overline{Char}_\chi\setminus Char_\chi$. Set $\Ohat(Char)^f=\bigoplus\limits
_{\chi} \Ohat(Char_\chi)$. Then it is easy to see that $\Ohat(Char)^f\cong \SS_b(\mT)$.
Choosing a point in $Char_\chi$ we get  identifications $Char_\chi\cong \mT^\vee$,
$\overline{Char}_\chi\cong \mTb^\vee$, 
$$\Ohat_\chi\cong \{ \phi: X_*(\mT)\to \Ce \ |\ \supp(\phi)
\subset S +X_*(\mT)^+\ {\mathrm{for}} \ S\subset X_*(\mT),\ |S|<\infty\}.$$
We define a linear functional $\int :\Ohat_\chi \to \Ce$ sending a function $\phi:X_*(\mT)\to \Ce$ to $\phi(0)$. It is easy to see that this is independent of the choice of a point in $Char_\chi$.
If $f\in \Ohat(Char_\chi)$ is a Laurent expansion of a rational function on $Char_\chi$
(we assume that the divisor of poles of $f$ does not pass through the zero dimensional 
orbit $\{t_0^\chi\}$ of $\mT^\vee$
on $\overline{Char}_\chi$, so that the Laurent expansion is well defined), then $\int f$
is the integral of $f$ over a coset of the maximal compact subtorus close to $t_0^\chi$, hence the notation.
We let $\int:\Ohat (Char)^f \to \Ce$ be the linear functional coinciding with the above functional
on $Char_\chi$ for every $\chi$.
It is easy to see that under the isomorphism $\Ohat(Char)^f\cong \SS_b(\mT)$, the functional $\int$
goes to the functional $h\mapsto h(1)$.

For a character $\chi$ of $T$ let $i_\chi:\H\to End(V_\chi)$, $i_\chi^-:\H\to End(V_\chi^-)$ denote the 
induced representations $i_B^G(\chi)$, $i_B^{-G}(\chi)$ respectively.

Recall from  \S \ref{61} that elements in $\MM(Y)$ can be interpreted as operators acting on $\SS(G/U)$ which commute with the action of the abstract Cartan $\mT$ and are invariant with respect to an open compact subgroup, or as family of operators acting on the space $V_\chi$  for every $\chi$. In particular, for $f\in \MM(Y)$ we get a regular function $\tau(f)\in \OO^f(Char)$
whose value at $\chi$ equals $Tr(f,V_\chi)$.

It is easy to see that $(\MM_b^1)_{sm}(Y)=\Ohat(Char)\otimes _{\O(Char)} \MM(Y)$, where
$(\MM_b^1)_{sm}\subset \MM_b^1$ denotes the subspace of elements invariant under some open subgroup in $G\times G$.
Thus
for $f\in (\MM_b^1)_{sm}(Y)$ we get $\tau(f)\in \Ohat(Char)$.

In the next statement we use Proposition  \ref{Iisom}.
 
\begin{Claim}\label{Pl}
a) For $f\in (\MM_b^1)_{sm}(Y)$ we have:
\begin{equation}\label{inttau}
\int \tau(f)= \int\limits _{\Delta_\BB\times \{1\}}f=\A^\star(f)(1),
\end{equation}
where $\Delta_{\BB\times \{1\} }\subset pr_Y^{-1}(\BB^2)\subset Y$ is 
defined by means of \eqref{pXo}.\footnote{Notice \label{foot4} that for 
$f\in (\MM_b^1)_{sm}(Y)$ the restriction $f|_{\Delta_\BB\times \{t\}}$,
$t\in \mT$
is a locally constant 
{\em measure} on $\BB$, thus the second integral in the displayed formula is well defined.}

b) For any $\phi\in \SS(X)$ 
we have $$\Be (\phi)(1) = \int \tau((I_{1}^\star)^{-1}(\phi)).$$ 
\end{Claim}

\proof Pick an open compact subgroup $K\subset G$ such that
$f\in (\MM_b^1)^{K\times K}(Y)$. The $\O(Char)$-module $\SS(G/U)^K$ splits
as a sum where each summand is free over the ring functions on a component of $Char$;
moreover, each module has a basis indexed by the finite set $K\backslash G/B$. 
For $f\in \MM(Y)$ the sum of diagonal entries of the corresponding operator on $\SS(G/U)$ in this basis is readily seen to equal the element
$$pr_{2*}(f|_{pr_Y^{-1}}(\Delta_\BB))\in \SS(\mT)\cong \O^f(Char),$$
 where $pr_2$ denotes projection to the second factor in \eqref{pXo} and
 we use  the observation in  footnote \ref{foot4} to make sense of $pr_{2*}$. Applying the functional 
 $\int:\O^f(Char) \to \Ce$ to this sum of diagonal elements yields $\int \tau(f)$, which is thus
 seen to be equal to the second expression in \eqref{inttau}. 
 This proves the first equality for $f\in \MM(Y)$, the case of $f\in (\MM_b^1)_{sm}(Y)$ follows.
 The second equality in (a)
  is clear since the preimage of $1\in G$ under the projection $(G\times G)/B\to G$
 maps isomorphically to $\BB\times \{1\}\subset pr^{-1}_Y(\Delta_\BB)\subset Y$. 
 
 This proves part (a); part  (b) follows 
from (a) compared to Corollary \ref{Bform}. \qed

\begin{Rem}
A similar equality holds where $I_{1}$ is replaced by $I_w$ for some $w\in W$. 
\end{Rem}

\subsection{Commutative subalgebras in the Hecke algebra}
We fix a maximal split torus $T\subset G$ and a subgroup $K\subset K_0\subset G$
(where $K_0$ is as in section \ref{53}) which is nice in the sense of \cite{BD} and is in a good relative position to $T$. Recall that this means that
\begin{equation}\label{KKK}
K=K^+\cdot K^0\cdot K^-
\end{equation}
 for any pair of opposite Borel subgroups $B^+=U^+T$, $B^-=U^-T$, where $K^+=U^+\cap K$, $K^-=U^-\cap K$, $K^0=T\cap K$. We fix a Borel subgroup containing $T$, this defines the cone of dominant coweights
$X_*(T)^+\subset X_*(T)$.

The following technical statement is not hard but it plays an important role in Bernstein's theory.

Recall  notation $\iota$, $\iota_I$ introduced in section \ref{53}.
For $\nu\in 
X_*(T)^+$ let  $e_\nu\in \H$ be the delta-function of the two-sided coset $K\iota(\nu)K$.

 Let $x_\nu\in \SS(X)$ be the delta-function of the $K\times K$ orbit of $\iota_X(\nu)$.
We set
$\theta_\nu=\Be (x_\nu)$.

\begin{Prop}\label{thetaprop}
 The map $X_*(T)\to \H$, $\nu\mapsto \theta_\nu$ is uniquely characterized by
the following two properties.

i) $\theta_{\mu_1+\mu_2}=\theta_{\mu_1} \theta_{\mu_2}$.

ii) There exists $\nu_0$ such that for $\nu \in \nu_0+X_*(T)^+$, we have $\theta_\nu =e_\nu$.
\end{Prop}

We will need the following Lemma. 

\begin{Lem}\label{munulem} \cite{BD}
a) For $\mu, \nu \in X_*(T)^+$ we have $e_\mu e_\nu=e_{\mu+\nu}$.

b) Given $\lambda \in X_*(T)^+$, there exists $N>0$ such that the kernel of left multiplication
by $e_{N\lambda}$ on $\H_K$ equals the kernel of left multiplication by $e_{(N+1)\lambda}$.
\end{Lem}

\proof (a) follows directly from \eqref{KKK}. 
 (b) follows from the fact that $\H_K$ is right Noetherian. \qed

\noindent{\bf{Proof}} of Proposition. To check uniqueness, assume  $\theta_\mu$ and $\theta_\mu'$ are different collections of elements
satisfying (i,ii). Fix $\nu\in X_*(T)^+$ such that $\theta_\nu=e_\nu=\theta'_\nu$. We have also
$\theta_{\mu+N\nu}=e_{\mu+N\nu}=\theta'_{\mu+N\nu}$ for some $N$. Thus $\theta_\mu-\theta'_\mu$ lies
in the image of left multiplication by $e_\nu$, as well as in the kernel of left multiplication
by $e_{N\nu}$. Lemma \ref{munulem}(b) implies that $\theta_\mu=\theta_\mu'$, which proves
uniqueness.

It remains to check that $\theta_\nu=\Be (x_\nu)$ does satisfy (i,ii). 
It follows from Lemma
\ref{cosets}(b) that given
 an admissible bijection between $U\subset G$ and $V\subset X$ there
exists $\nu_0\in X_*(T)^+$ such that the bijection sends the set $K\iota_G(\nu)K$
into the set $K\iota_X(\nu) K$ for $\nu \in \nu_0+X_*(T)^+$. Thus (ii) follows from the construction of $\Be $. Given $\nu\in X_*(T)^+$, Proposition \ref{compatibility} implies that
there exists $\nu_0\in X_*(T)^+$ such that
\begin{equation}\label{exe}
e_\nu*x_\mu=x_{\nu+\mu}=x_\mu*e_\nu
\end{equation}
 for $\mu \in \nu_0+
X_*(T)^+$. Since the elements $x_\mu$ are permuted by the action of $\mT$ commuting with the action of
$\H\otimes \H$, we see that \eqref{exe} holds for all $\mu$.
Hence $e_\nu\theta_\mu=\theta_{\nu+\mu}$ for all $\mu$. We now proceed to check (i).
Pick $\nu$ such that $e_\eta=\theta_\eta$ when $\eta$ is either 
 $\nu+\mu_1$, or $\nu+\mu_1+\mu_2$. Then we see that $e_\nu\theta_{\mu_1+\mu_2}=e_{\nu+\mu_1+\mu_2}=e_\nu\theta_{\mu_1}
\theta_{\mu_2}.$ Also, both $\theta_{\mu_1+\mu_2}$ and $\theta_{\mu_1}$, $\theta_{\mu_2}$
lie in the image of left multiplication by $e_{N\nu}$ for all $N$. Thus the desired equality follows
from Lemma \ref{munulem}(b). \qed

\begin{Ex}\label{thetaex} The case most often encountered in the literature is when $K=I$ is the Iwahori subgroup.
Then $\H=\H_I$ is isomorphic to the affine Hecke algebra, it contains invertible elements
$\theta_\nu^I$, $\nu\in \Lambda$ such that $\theta_{\nu_1+\nu_2}^I=\theta_{\nu_1}\cdot \theta_{\nu_2}^I$
and $\theta_\nu^I=\delta_{I\iota(\nu) I }$ when $\nu$ is dominant;
they form
a part of a system of generators for this algebra discovered by Bernstein \cite{Lus}, which have
found numerous important applications. It is clear that these are exactly the elements described
in Proposition \ref{thetaprop} in the case $K=I$.

We point out a special case of this potentially useful in applications.
Setting $\nu=0$ we get:
$$\Be (\delta_{I\iota_X(0)I})=\theta_0^I=\delta_I.$$ 
Furthermore, applying averaging with respect to $K_0\times K_0$ to both sides and using that
$X(O)=K_0\cdot \iota_X(0)\cdot K_0$ we get:
$$\Be (\delta_{X(O)})=\delta_{K_0}.$$

\end{Ex}


\subsection{Proof of Proposition \ref{BstimDr}}
It is easy to see that 
$$X\cap \Xb_{\aff}(O)=\bigcup \limits_{\la\preceq_\Qu 0} X_\la,$$
where we write $\la\preceq_\Qu \mu$ if $\mu-\lambda$
is a linear combination 
of positive coroots with nonnegative rational coefficients. 

Using the notation and the content  of Example \ref{thetaex} we reduce  the claim to showing that 
$$\supp(\theta_\lambda^I)\cap K_0=\emptyset\ \ \ \ \ {\mathrm unless}\ \ \ \ \ \lambda\preceq_\Qu 0.$$

In fact, we will show that
\begin{equation}\label{suppt}
\supp(\theta_\lambda^I)\cap K_0=\emptyset\ \ \ \ \ {\mathrm unless}\ \ \ \ \ \lambda\preceq 0,
\end{equation}
where $\preceq$ is
 standard order on weights, i.e. $\lambda\preceq \mu $
if $\mu-\la$ is a sum of positive roots.

We have the standard basis $t_w$ of $\H_I$, $w\in W_{\aff}$ where $W_{\aff}$ is the extended affine Weyl group. We also have an algebra antiinvolution 
$i:\H_I\to \H_I$ sending $t_w$ to $t_w^{-1}$ and 
a pairing $\H_I\times \H_I\to k$ given by $\langle f,g\rangle =
c_1$ where $i(f)g=\sum c_w t_w$. Then \eqref{suppt} reduces to:
\begin{equation}\label{suppt1}
\langle  t_w, \theta_\lambda^I \rangle =0  \ \ \ \ {\mathrm for}\ \ \ \ \ w\in W,\ \ \ \lambda\not \preceq 0,
\end{equation}
where $W$ is the finite Weyl group.

Pick a dominant weight $\mu$ such that $\mu +\lambda$ is also dominant. 
Then for $w\in W$ we have:
$$\langle t_w, \theta_\lambda^I\rangle = \langle t_\mu t_w, t_\mu \theta_{\lambda}\rangle
=\langle  t_{(\mu)\cdot w}, \theta_{\mu+\lambda}^I\rangle.$$
It is easy to see that for $w,v\in W_{\aff}$
we have $\langle t_w,t_v\rangle=0$ unless $v\leq w$ where 
$\leq$ refers to the Bruhat order on $W_{\aff}$.

For dominant coweights  $\eta$, $\nu$ the inequality $\eta \leq (\nu)w$, $w\in W$ implies
that $\eta\preceq \nu$. So$ \langle  t_w, \theta_\lambda^I \rangle \ne 0$ implies
that $\mu+\lambda\preceq \mu$, i.e. $\lambda \preceq 0$, which proves \eqref{suppt1}
and hence Proposition \ref{BstimDr} for adjoint groups. \qed

\subsection{Plancherel functional on the abelian subalgebra}
 Corollary \ref{Bform} shows that
\begin{equation}\label{AI1}
\theta_\nu(g)=\Be(x_\nu)(g)=\A^\st I_1 (g).
\end{equation}
This  yields a spectral expression for the Plancherel functional $f\mapsto f(1)$ restricted to the subalgebra
$\A_K\subset \HH$ spanned by $\theta_\nu$.
To state it we introduce the following notation.


We define a rank one operator $\Pi:V_\chi^K\to (V_\chi^-)^K$ as follows.
Notice that the space $V_\chi^K$ splits as a direct sum of one dimensional subspaces
indexed by the set $K\bs G/B$; likewise, $V^-_\chi$ is a direct sum of one dimensional subspaces
indexed by $K\bs G/B^-$. The one dimensional summands corresponding to the coset of 1 in the
two spaces are canonically isomorphic (both are identified with the space of the one dimensional
representation $\chi$ of $T$), let us denote this space by $\Ce_\chi$.
 We define $\Pi:(V^-_\chi) ^K\to V_\chi^K$ as the composition
$(V^-_\chi)^K\to \Ce_\chi \to V_\chi^K$, where the first arrow is the projection arising
from the above splitting into a direct sum, and the second arrow is the embedding.


In the next Theorem we consider a function $\chi\mapsto Tr(h \Pi R^{-1}, V_\chi)$, $h\in \HH$,
where $R:V^-_\chi \to V_\chi$ is the intertwining operator defined for $\chi$ in a Zariski open dense
subset of $Char$.
This is a rational function on $Char$ vanishing on all but a finite number of components. Moreover,
in view of Proposition \ref{sigmasi}, for every $\chi$ the point $t_0\in \overline{Char}_\chi$ does not lie
in the divisor of poles of this function. Thus the functional $\int$ introduced in
\S \ref{101} is well defined on such a function.

\begin{Thm}\label{flaThm}
For $h\in A_K$ we have
$$h(1)=\int Tr(h \Pi R^{-1}, V_\chi^K).$$
\end{Thm}

\proof We can assume that $h=\theta_\nu$.
Applying Corollary \ref{Bform} with $w=1$, we get $\theta_\nu=\A^* I^{-1}(x_\nu)$.

A function $\phi\in \MM_b^{w=1}(Y)$ defines an operator on the completion of the universal principal series, $\SS(G/U)\otimes _{\OO(Char)}\Ohat(Char)$.

By Claim  \ref{Pl}(a) we have 
$$\int Tr(\phi, V_\chi) =A^*(\phi)(1),$$
where  integral in the left hand side
has the same meaning as in the statement of the Theorem. 

Thus the proof will be complete once we show that
\begin{equation}\label{iIx}
i_\chi(I_{1}(x_\nu))=i_\chi(\theta_\nu) \Pi  R^{-1}.
\end{equation}

A function $f\in \SS(X)$ defines an operator between the principal series and the opposite principal series, $\sigma_\chi(f):V^-_\chi \to V_\chi$. In particular, it is easy to see that $\sigma_\chi (x_0)=\Pi$ for all $\chi$.
Also, it follows from the definitions that $i_\chi((I_{1}^\star)^{-1}(f))=R^{-1}\circ \sigma_\chi(f)^{-1}$.
Finally, using that all the maps involved commute with the $\H\otimes \H$ action, we get:
$$i_\chi((I_{1}^\star)^{-1}(x_\nu))=i_\chi(\theta_\nu*x_0)=i_\chi(\theta_\nu) \Pi R^{-1}.$$
Thus we have proven \eqref{iIx}, and hence the Theorem. \qed

\begin{Rem}
  In the special case  when $K=I$ is the Iwahori
subgroup  Theorem \ref{flaThm} reduces to (an equivalent form) of 
  \cite[Theorem 1.14]{Opdam1}. 
\end{Rem}

\begin{Rem}
We expect that a version of Claim \ref{Pl}(b) holds for the space of rapidly decreasing functions on $X$ and that it can be used to deduce the standard   Plancherel formula. We plan to develop this
application in the future.
\end{Rem}
\section{Appendix: quasi-normal cone for toric coverings}

\bigskip

\centerline{by ROMAN BEZRUKAVNIKOV, DAVID KAZHDAN}
\centerline{ and YAKOV VARSHAVSKY}


\bigskip

In this section we introduce a version of the definition of a normal cone which behaves well for a class
of singular varieties, including De Concini - Procesi compactifications of not necessarily adjoint groups.

\subsection{Quasi-normal cone}\label{qnc} Let $X$ be a smooth variety over a field $F$ and $Z\subset X$ be a smooth locally closed subvariety. Let $X'$ be a normal scheme and $X'\to X$ be a finite flat morphism, let $Z'$ be the preimage of $Z$ equipped with reduced subscheme structure. Let $N_X(Z)$ be the normal bundle.

Let $\tN_X(Z)$ be the deformation  to the normal cone; thus
$(\Aone \setminus \{0\})\times X\subset \tilde N_X(Z)\supset
N_X(Z)$. Recall that $\Nt$ comes equipped with a $\Gm$ action which dilates the fibers of the normal
bundle and acts on $(\Aone \setminus \{0\})\times X$ by $t:(z,x)\mapsto (tz,x)$.

Let $\tN_{X'}(Z')$ be the normalization of $\tN_X(Z)$ in $(\Aone \setminus \{0\})\times X'$, and set $N_{X'}(Z')=\tN_{X'}(Z)\times_{\Aone}\{0\}$.
We call $N_{X'}Z'$ the {\em quasi-normal cone} of $Z'$ in $X'$.

We have a locally closed embedding $Z'\times \Aone \to \Nt_{X'}Z'$.

{\bf Definition.} We say that $X'$ is {\em well approximated around $Z$} if the following conditions
hold.

\begin{enumerate}\item
 For some $d>0$ the composition of the natural action of $\Gm$ on $\tN_Z(X)$ with the map
$\Gm \to \Gm$, $t\mapsto t^d$, lifts to an action $\alpha_d$ of $\Gm$ on $\Nt_{Z'}(X')$.

\item Zariski locally on $Z$ there exists an isomorphism $\tau$ between the formal neighborhoods of
$Z'\times \Aone$ in $N_{X'}(Z')\times \Aone$ and in $\Nt_{Z'}(X')$, such that

\begin{enumerate}
\item $\tau$ restricted to the preimage of $0\in \Aone$ equals identity.

\item $\tau$ is $\Gm$-equivariant where
$\Gm$ acts  on $\Nt_{Z'}(X')$  via $\alpha_d$ and on
 $N_{X'}(Z')\times \Aone$ by $t:(x,z)\mapsto (\alpha_d(x),t^dz)$.
\end{enumerate}
\end{enumerate}

\medskip

If $F$ is a local field, we can repeat the definition replacing an isomorphism of formal neighborhoods
by an analytic isomorphism of actual $F^\times$ invariant neighborhoods in the space of $F$-points;
 "Zariski local"
should then be replaced by local in the sense of $F$-topology. We will use the term "analytically
well approximated" for this version of the property.
If $F$ is non-Archimedian, existence of local isomorphisms implies existence of a global isomorphism of appropriate neighborhoods.

In the latter case, the restriction of $\tau$ to the fiber over $1\in \Aone$ will be called an {\em admissible bijection}.

\begin{example}
If $X'=X$, then $\tau$  amounts to  an isomorphism between a neighborhood of $Z$ in $X$ and a neighborhood
of the zero section in the normal bundle $N_X(Z)$, whose differential in the normal direction equals identity.
\end{example}

\subsection{Toric covering}
Let $X$ be again a smooth variety over $F$, and let $ D_i, i\in I$ be smooth divisors with normal
crossing in $X$. For a subset $J\subset I$
we have the corresponding stratum  $X_J = \cap D_i, i\in J.$ We fix $J\subset I$ and let $Z=X_J$.

We have line bundles $\LL_i = \OO(D_i)$ on $X$, each coming with a canonical section $s_i$. They can be combined into a $T$-bundle $\EE$ for the torus $T = (\Gm)^I.$

Let $T'$ be another torus and $T'\to T$ be a fixed isogeny. Suppose that
the above $T$-bundle on $X$ lifts to  a $T'$-bundle $\EE'$ (i.e. $\EE$ is the push-forward of $\EE'$); we fix such a $T'$-bundle $\EE'$. This data defines a ramified covering $X'\to X$ as follows.

Let $A=(\Aone)^I$, thus $T$ acts on $A$ making it a toric variety.  The $T$-bundle $\EE$ defines an associated bundle $A_\EE$ with the fiber $A$;
of course $A_\EE$ is nothing but the total space of the sum of line bundles $\LL_i$; the sections $s_i$
combine to a section $s:X\to A_\EE$.

Let $A'$ be the normalization of $A$
in the covering $T'\to T$. We also can form an associated bundle $A'_{\EE'}$.
Set $X'=X\times _{A_\EE} A'_{\EE'}$. We call $X'$ obtained by this construction a toric covering.

For a subset $J\subset I$ we let $A'(J)$ be the normalization of $(\Aone)^J$ in the quotient
of $T'$ by the neutral connected component  of the preimage of $(\Gm)^{I\setminus J}\subset T$,
let $A'_{\EE'}(J)$ be the corresponding  bundle over $X$.

\begin{example}\label{exl1} Let $G$ be a 
 reductive algebraic group.
  Let $G_{ad}=G/Z(G)$, where $Z(G)\subset G$ is the center, and let $G'\subset G$ be the derived
(commutator) subgroup.
Let $X=(G/G')\times \Gb_{ad}$ where $\Gb_{ad}$ is the wonderful compactification of $G_{ad}$, and $D_i$ be the components of the complement
$X\setminus (G/G')\times G_{ad}$. The torus $T$ is identified with the (abstract) Cartan subgroup of $G_{ad}$.
The total space of $\EE$ is easily seen to be identified with $(G/G')\times S^0(G_{ad})$, where
$S^0(G_{ad})\subset S(G_{ad})$ is an open part in Vinberg's semi-group $S(G_{ad})$.
 Let $T'$ be the abstract Cartan of $G'$, and let $\EE$ be the $T$-bundle whose total space is $S^0(G)$. [If $G$ is simply-connected, this is a universal torus bundle
over $X$, i.e. the map from $X^*(T)$ to $Pic(X)$ given by this bundle is an isomorphism; we neither
 prove nor use this fact].
It is easy to see that the preimage of $(G/G')\times G_{ad}\subset X$ in $X'$ is isomorphic to $G$.
It follows from the Theorem below that $X'$ is the normalization of $(G/G')\times \Gb_{ad}$ in $G$.
\end{example}

 \begin{example}\label{exl2}
 Assume that the line bundles $\LL_i$ are trivial (i.e. $D_i$ is cut out by a global function). Fix $\EE'$ to be the trivial $T'$-torsor. Then $s$ amounts to a smooth morphism
$X\to A$, so that $Z$ is the preimage of the closure of a $T$-orbit in $A$; we have $X'=X\times_A A'$.
\end{example}

\subsection{Theorem}\label{qnrml}
{\em Let $X'\to X\supset D_i$ be a toric covering and $Z=X_J$ as above.

a) $X'$ is a normal variety and the map $X'\to X$ is finite and flat.

b) $X'$ is well approximated around $Z$. If $F$ is a local field it is also analytically well approximated.

c) The quasi-normal cone $N_{X'}(Z')$ is canonically isomorphic to
$A'_{\EE'}(J)\times _X Z$.}

\medskip

{\bf Remark.} If $T'=T$, so that $X'=X$, the isomorphism of part (c) amounts to the adjunction formula
$\OO(D)|_D\cong N_X(D)$.

\medskip

{\bf Proof.} 
It is easy to check that all the definitions above are compatible with smooth base change.
Thus it suffices to check that the statements hold after base change with respect to the morphism
from the total space of the torsor $\EE$ to $X$. This reduces the proof to the situation
of Example \ref{exl2}. Using compatibility with smooth base change again, we see that it suffices
to show that $A'$ is normal, $A'\to A$ is flat and finite; that $A'$ is well approximated around
the closure $A_J$ of a $T$-orbit in $A$, and that $N_{A'}(A_J')= A'$.
Normality and finiteness are clear from the definition, and flatness is easy to show.
The deformation to the normal cone $\tN_A(A_J)$ is constant, which implies the rest. \qed

\medskip

{\bf Remark.} Part (c) of the Theorem shows that an open part of $N_{X'}(Z')$ is identified
with the total space of $\EE'|_Z$. In the situation of Example \ref{exl1}, this is the space denoted
by $X_J$ in the paper.

\footnotesize{ {\bf R.B.}: Department of Mathematics, Massachusetts Institute
of Technology,\\ Cambridge MA 02139, USA; 
}

\footnotesize{ {\bf R.B.}:  
National Research University Higher School of Economics,
International Laboratory of Representation
Theory and Mathematical Physics,
20 Myasnitskaya st., Moscow 101000, Russia}


\footnotesize{ {\bf D.K.}: Department of Mathematics, Hebrew University, Jerusalem, Israel
}

\footnotesize{ {\bf Y.V.}: Department of Mathematics, Hebrew University, Jerusalem, Israel
}

\end{document}